\documentclass[12pt]{smfart}
\usepackage{amsmath,amscd,latexsym,verbatim,amssymb}
\usepackage[frenchb]{babel}     
\usepackage[all]{xy}
\usepackage{times}

\renewcommand{\hom}{{\operatorname{hom}}}        
\newcommand{\num}{{\operatorname{num}}} 

\newcommand{\kim}{{\operatorname{kim}}}

\newcommand{\gm}{{\operatorname{gm}}}
\newcommand{\Spec}{\operatorname{Spec}}
\newcommand{\Pic}{\operatorname{Pic}}

\newcommand{\ord}{\operatorname{ord}}
         
\newcommand{\resp}{{\it resp.} }
\newcommand{\rg}{\operatorname{rg}}
\newcommand{\Cores}{\operatorname{Cor}}
\newcommand{\Chow}{\operatorname{Chow}}
\newcommand{\Tot}{\operatorname{Tot}}
\newcommand{\cont}{{\operatorname{cont}}}

\newcommand{\el}{{\operatorname{ell}}}
\newcommand{\ab}{{\operatorname{ab}}}
\newcommand{\eff}{{\operatorname{eff}}}
\newcommand{\ind}{{\operatorname{ind}}}
\newcommand{\tate}{{\operatorname{tate}}}
\newcommand{\tors}{{\operatorname{tors}}}
\newcommand{\cotors}{{\operatorname{cotors}}}

\newcommand{\oo}{\operatornamewithlimits{\otimes}}
\newcommand{\cf}{{\it cf.} }

\newcommand{\loccit}{{\it loc. cit.} }
\newcommand{\Zar}{{\operatorname{Zar}}}
\newcommand{\et}{{\operatorname{\acute{e}t}}}
\newcommand{\Et}{{\operatorname{\acute{E}t}}}
\newcommand{\If}{\Longrightarrow}

\newcommand{\G}{\mathbb{G}}
\newcommand{\N}{\mathbf{N}}
\renewcommand{\P}{\mathbf{P}}
\newcommand{\Q}{\mathbf{Q}}

\newcommand{\Z}{\mathbf{Z}}

\newcommand{\sA}{{\mathcal{A}}}

\newcommand{\sC}{{\mathcal{C}}}
\newcommand{\sD}{{\mathcal{D}}}
\newcommand{\sF}{{\mathcal{F}}}

\newcommand{\sH}{{\mathcal{H}}}

\newcommand{\sK}{{\mathcal{K}}}
\newcommand{\sM}{{\mathcal{M}}}

\newcommand{\sO}{{\mathcal{O}}}

\newcommand{\sT}{{\mathcal{T}}}

\newcommand{\sV}{{\mathcal{V}}}

\newcommand{\To}{\longrightarrow}

\newcommand{\Inj}{\lhook\joinrel\longrightarrow}
\newcommand{\surj}{\rightarrow\!\!\!\!\!\rightarrow}
\newcommand{\Surj}{\relbar\joinrel\surj}
\newcommand{\iso}{\overset{\sim}{\longrightarrow}}
\newcommand{\osi}{\overset{\sim}{\longleftarrow}}

\newcommand{\un}{{\bf 1}}

\renewcommand{\epsilon}{\varepsilon}
\renewcommand{\phi}{\varphi}

\newcommand{\colim}{\operatornamewithlimits{\varinjlim}}
\newcommand{\plim}{\operatornamewithlimits{\varprojlim}}

\newcommand{\Ker}{\operatorname{Ker}}

\newcommand{\Coker}{\operatorname{Coker}}

\newcommand{\car}{\operatorname{car}}

\newcounter{spec}
\newenvironment{thlist}{\begin{list}{\rm{(\roman{spec})}}%
{\usecounter{spec}\labelwidth=20pt\itemindent=0pt\labelsep=10pt}}%
{\end{list}}

\newtheorem{Thm}{Th\'eor\`eme}
\newtheorem{Cor}{Corollaire}
\theoremstyle{definition}
\newtheorem{Defn}{D\'efinition}
\newtheorem{Ex}{Exemples}

\swapnumbers

\theoremstyle{theorem}
\newtheorem{thm}{Th\'eor\`eme}[section]
\newtheorem{lemme}[thm]{Lemme}
\newtheorem{prop}[thm]{Proposition}

\newtheorem{cor}[thm]{Corollaire}

\theoremstyle{definition}
\newtheorem{defn}[thm]{D\'efinition}

\newtheorem{qn}[thm]{Question}

\newtheorem{rem}[thm]{Remarque}
\newtheorem{rems}[thm]{Remarques}

\numberwithin{equation}{section}

 at10pt

\newcommand{\prf}{\noindent {\bf D\'emonstration. }}
\renewcommand{\qed}{\hfill $\square$\medskip}

\begin{document}

\title[Produits de courbes elliptiques sur un corps fini]
{\'Equivalence rationnelle, \'equivalence num\'erique et produits de
courbes elliptiques sur un corps fini}
\author{Bruno Kahn}
\address{Institut de Math\'ematiques de
Jussieu\\175--179 rue du
Chevaleret\\ \break 75013
Paris\\France.}
\email{kahn@math.jussieu.fr}
\date{20 juillet 2002}

\renewcommand{\abstractname}{Abstract}

\begin{abstract}
We prove that if $X$ is a product of elliptic curves over a
finite field $k$, rational and numerical equivalences agree on $X$. The
proof works for any smooth projective variety $X$ ``of abelian type" for
which the Tate conjecture holds (the case of products of elliptic curves
being due to M. Spie\ss's \protect\cite{spiess}). It uses Soul\'e's ideas
\protect\cite{soule}, U. Jannsen's semi-simplicity theorem
\protect\cite{jannsen},   and a result of Y. Andr\'e and the
author \protect\cite{ak} inspired by S.I. Kimura's results on
finite-dimensional Chow motives \protect\cite{ki}. We give some
consequences, among which: the conjectures of Lichtenbaum
\protect\cite[\S 7]{licht} hold true for $X$,
the second Chow group of $X$ is finitely generated, the Beilinson-Soul\'e
conjecture holds in weight $n$ for the function field of $X$ provided
$n\le 2$ or $\dim X\le 2$, Gersten's conjecture holds for discrete
valuation rings with residue field such a function field if $\dim X\le 2$.
Also, if
$U$ is an open subset of
$X$ with
$\dim X\le 2$, the action of Frobenius on $H^*_c(\bar U,\Q_l)$ is
semi-simple and its characteristic polynomial does not depend on $l\ne
\text{ char } k$.
\end{abstract}
\maketitle

\enlargethispage*{20pt}

\section*{Introduction}

Soient $k$ un corps fini de caract\'eristique $p$ et $G=Gal(\bar k/k)$.
Dans un article s\'eminal \cite[th. 4 i)]{soule}, Soul\'e
d\'emontre que pour toute vari\'et\'e projective lisse $X$ de dimension
$\le 3$ et ``de type ab\'elien" (par exemple une vari\'et\'e ab\'elienne
ou un produit de courbes), l'application ``classe de cycle $l$-adique"
\begin{equation*}\label{eq0}\tag{*}
CH^i(X)\otimes\Q_l\to H^{2i}(\bar X,\Q_l(i))^G 
\end{equation*}
est bijective pour tout nombre premier $l\ne p$ et tout
$i\ge 0$. Ceci est un cas particulier de deux conjectures fondamentales
en g\'eom\'etrie arithm\'etique: la conjecture de Tate \cite{tate3}
(surjectivit\'e de
\eqref{eq0}) et la conjecture de Beilinson
\cite[1.0]{beilinson}\footnote{Mais aussi de Lichtenbaum en conjonction
avec la conjecture de Tate, \cf\protect\cite[\S 7, 7)]{licht}.}
(injectivit\'e de
\eqref{eq0}). La d\'emonstration consiste \`a consid\'erer la
d\'ecomposition motivique de
$X$ pour se r\'eduire au th\'eor\`eme de Tate sur les endomorphismes de
vari\'et\'es ab\'eliennes
\cite{tate2} et \`a l'hypoth\`ese de Riemann sur $k$
(Deligne \cite{deligne}).

La conjecture de Tate est maintenant d\'emontr\'ee pour un bon nombre
de vari\'et\'es ab\'eliennes sur $k$ (voir ci-dessous). La situation est
diff\'erente pour la conjecture de Beilinson, qui n'est connue
grosso-modo qu'en codimension $1$ (o\`u elle est triviale); elle est
d'ailleurs utilis\'ee par Soul\'e dans ce cas pour obtenir \eqref{eq0}.
Dans le pr\'esent article, nous uti\-li\-sons une id\'ee nouvelle de S.I.
Kimura
\cite{ki} pour la d\'emontrer en toute codimension dans un certain
nombre de cas nouveaux. 

Plus pr\'ecis\'ement, notant $\sA(k)=\sA$ la cat\'egorie des motifs de
Chow sur $k$ \`a coefficients rationnels:

\begin{Defn} a) Soit $B(k)$ l'ensemble des classes d'isomorphismes de
$k$-vari\'et\'es projectives lisses dont le motif de Chow est dans la
sous-cat\'egorie
\'epaisse rigide de $\sA$ engendr\'ee par les motifs d'Artin et les
motifs de vari\'et\'es ab\'eliennes (ou de courbes, c'est la m\^eme
chose).\\
b) Soit $B_\tate(k)\subset B(k)$ le sous-ensemble des vari\'et\'es
$X$ v\'erifiant la conjecture de Tate (pour un nombre premier $l\ne p$
donn\'e: cela ne d\'epend pas de $l$, \cf \cite[Th. 2.9]{tate2} puisque
Frobenius op\`ere de mani\`ere semi-simple sur la cohomologie
$l$-adique de $X$).
\end{Defn}

\begin{Ex} a) $B(k)$ contient l'ensemble $A(k)$  d\'efini dans
\cite[3.3.1]{soule}.\\
b)  $B_\tate(k)$ contient les $X\in B(k)$ avec $\dim
X\le 3$ (Soul\'e, voir ci-dessus) et les produits de courbes elliptiques
\cite{spiess}.\\
c) Milne a d\'emontr\'e que la conjecture de Hodge pour les
vari\'et\'es ab\'eliennes complexes de type $CM$ implique la
conjecture de Tate pour les vari\'et\'es ab\'e\-lien\-nes sur $k$
\cite{milne5}. Inconditionnellement, les r\'esultats de
\cite{milne4} montrent qu'on peut prouver la conjecture de
Tate pour ``beaucoup" de vari\'et\'es ab\'eliennes $A$ sur $k$. C'est le
cas si l'alg\`ebre des ``cycles de Tate" de $A$ est
engendr\'ee en degr\'e $1$, par r\'eduction au th\'eor\`eme de Tate
\cite{tate2}. C'est ainsi que Spie\ss\ d\'emontre la conjecture de Tate
pour les produits de courbes elliptiques. D'autres
exemples sont les puissances de vari\'et\'es ab\'eliennes simples ``de
type K3" de Zarhin \cite{zarhin} ou ``presque ordinaires" de
Lenstra-Zarhin \cite{len-zar}, \cf \cite[A.7, exemples]{milne4}. Pour un
exemple o\`u cette condition n'est pas v\'erifi\'ee mais o\`u n\'eanmoins
la conjecture de Tate est d\'emontr\'ee, voir \cite[ex.
1.8]{milne4}.\\
d) Un exemple amusant de nature tr\`es l\'eg\`erement diff\'erente est
celui d'une hypersurface de Fermat
\[X: X_0^m+\dots+X_{d+1}^m=0\]
o\`u $m$ est tel que $q^\nu\equiv -1\pmod{m}$ pour un $\nu$ convenable
(Tate, Katsura-Shioda, \cite{tate3,ks}).
\end{Ex}

\begin{Thm}[\cf th. \protect{\ref{t1}}]\label{T1} Pour tout
$X\in B_\tate(k)$, l'\'equivalence ra\-tion\-nel\-le est \'egale \`a
l'\'e\-qui\-va\-len\-ce num\'erique (\`a coefficients rationnels) sur $X$.
\end{Thm}

Dans \cite{tate}, nous avons montr\'e que la conjonction des
conjectures de Tate et de Beilinson a des implications
consid\'erables. Dans cet esprit, nous tirons les cons\'equences suivantes
du th\'eor\`eme \ref{T1}. Soient $X\in B_\tate(k)$, $d=\dim X$ et
$K=k(X)$.

\begin{sloppypar}
\begin{Cor}[cor. \protect{\ref{c3l}}]\label{T2} Les conjectures de
Lichtenbaum
\cite[\S 7]{licht} sont vraies pour $X$.
\end{Cor}
\end{sloppypar}

\begin{Cor}[th. \protect{\ref{c6}}]\label{T4} Pour tout ouvert
$U$ de $X$, les groupes $H^0(U,\sK_2)$, $H^1(U,\sK_2)$ et
$H^2(U,\sK_2)=CH^2(U)$ sont de type fini. Le groupe
$K_3(K)_\ind$ est \'egal \`a $K_3(k)$.
\end{Cor}

\begin{Cor}[cor. \protect{\ref{c2}}]\label{T3} La conjecture de
Beilinson-Soul\'e est vraie en poids $n$ pour $K$ pourvu que $n\le 2$ ou
que $d\le 2$.
\end{Cor}

\`A ma connaissance, ce sont les premiers exemples non triviaux (au-del\`a
de la dimension 1) o\`u cette conjecture est explicitement d\'emontr\'ee.
Toutefois, pour $d=2$ ou pour $n=2$, $d=3$, les r\'esultats de
\cite{soule} suffiraient (voir \cite[th. 4 iii)]{soule} et la
d\'emonstration du corollaire \ref{c2}).

\begin{Cor}[th. \protect{\ref{c4}}]\label{T5} Si $d=2$, on a des
isomorphismes canoniques pour tout $n\ge 0$
\[\left(K_n^M(K)\oplus\bigoplus_{0\le i\le n-1}
H^{2i-n-1}(K,(\Q/\Z)'(i))\right)\otimes\Z_{(2)}\iso
K_n(K)\otimes\Z_{(2)}\]  o\`u $(\Q/\Z)'(i)=\colim\limits_{(m,\car k)=1}
\mu_m^{\otimes i}$.
\end{Cor}

Si l'on fait ``tendre $k$ vers l'infini" dans le th\'eor\`eme \ref{T5},
on obtient une confirmation partielle d'une conjecture de Suslin
\cite[conj. 4.1 et note]{suslinicm}: si un corps $F$ contient un
corps alg\'ebriquement clos $F_0$, $K_*(F)$ est
engendr\'e multiplicativement par $K_1(F)$ et $K_*(F_0)$ (corollaire
\ref{c4.2} et remarque \ref{r2}). Ce r\'esultat est faux en
caract\'eristique z\'ero d'apr\`es de Jeu \cite{jeu}.

\begin{Cor}[cor. \protect{\ref{c4.1}}]\label{T5bis} Supposons $d=2$, et
soit $A$ un anneau de valuation discr\`ete de corps des fractions $E$ et
de corps r\'esiduel
$K$. Alors la conjecture de
Gersten est vraie pour la $K$-th\'eorie alg\'ebrique de $A$:  pour tout
$n\ge 0$, la suite
\[0\to K_{n+1}(A)\to K_{n+1}(E)\to K_n(K)\to 0\]
est exacte.
\end{Cor}

\begin{Cor}[prop. \protect{\ref{p6}}] \label{T6} Supposons $d=2$, et
soit $U$ un ouvert de $X$. Alors, pour tout $l\ne p$, l'action de
Frobenius sur
$H^*_{c}(\overline U,\Q_l)$ est semi-simple. De plus, le polyn\^ome
caract\'eristique de cette action est in\-d\'e\-pen\-dant de $l$.
\end{Cor}

Le corollaire \ref{T6} est conjectur\'e (Grothendieck, Serre) pour tout
$k$-sch\'ema de type fini $U$. Gabber m'en a indiqu\'e une
d\'emonstration directe. La d\'e\-mons\-tra\-tion donn\'ee ici consiste
\`a montrer que le motif \`a supports compacts de $U$ dans la cat\'egorie
triangul\'ee des motifs
$DM_\gm(k)\otimes\Q$ de Voevodsky \cite{voetri} est somme directe de
motifs purs d\'ecal\'es; elle peut servir de mod\`ele pour une
d\'e\-mons\-tra\-tion future du corollaire
\ref{T6} en g\'en\'eral.

Pour obtenir les corollaires \ref{T2} et \ref{T4}, nous
prouvons des r\'esultats plus pr\'ecis faisant intervenir une nouvelle
cohomologie introduite par Lichtenbaum \cite{licht2} (corollaire
\ref{c3}).
\bigskip

La d\'emonstration du th\'eor\`eme \ref{T1} utilise trois ingr\'edients
de mani\`ere essentielle: le fait que la conjecture de semi-simplicit\'e
est vraie pour la cohomologie $l$-adique de $X$, la semi-simplicit\'e de
la cat\'egorie des motifs modulo l'\'e\-qui\-va\-len\-ce num\'erique due
\`a  Jannsen \cite{jannsen} et un raffinement d'un r\'esultat de 
Kimura \cite[prop. 7.5]{ki} d\^u \`a Andr\'e et \`a l'auteur
\cite[prop. 9.1.14]{ak}. La technique de d\'emonstration, quant
\`a elle, remonte \`a Soul\'e \cite{soule} via Geisser \cite{geisser}.

Le corollaire \ref{T2} se d\'eduit du th\'eor\`eme \ref{T1} de
mani\`ere relativement classique \cite{milne2,milne3}; toutefois,
l'introduction de la cohomologie de Lichtenbaum mentionn\'ee ci-dessus
simplifie bien les choses. Il faut un peu d'effort pour d\'eduire la
version ``glo\-ba\-le" des conjectures de Lichtenbaum de leur version
localis\'ee en chaque nombre premier: \`a cet \'egard, le lemme
\ref{l3} est fort utile.

Le corollaire \ref{T3}, quant \`a lui, se d\'eduit assez facilement du
r\'esultat de Geisser selon lequel la conjugaison des conjectures de Tate
et de Beilinson implique la conjecture de Parshin (\cf corollaire
\ref{t2}).

La m\'ethode de d\'emonstration du th\'eor\`eme \ref{T1} et du corollaire
\ref{T2} conduit \`a de nouvelles formulations des trois conjectures dont
nous avons d\'emontr\'e l'\'equivalence dans \cite[th. 3.4]{glr}: voir
th\'eor\`eme \ref{t6}.

Cet article est construit comme suit. Au \S \ref{ratnum}, nous
d\'emontrons le th\'eor\`eme \ref{T1}. Dans les paragraphes suivants,
nous ``tirons les marrons du feu": cette op\'eration demande \`a
l'occasion le port de gants ignifuges. Comme indiqu\'e ci-dessus, les
con\-s\'e\-quen\-ces que nous donnons sont toutes des cas particuliers de
con\-s\'e\-quen\-ces des trois conjectures \'equivalentes de \cite[th.
3.4]{glr}, dont beaucoup ont \'et\'e d\'ej\`a d\'egag\'ees dans
\cite{tate}: nous nous sommes efforc\'e de r\'ediger les d\'emonstrations
de telle mani\`ere qu'elles puissent \^etre reproduites telles quelles
quand ces conjectures seront d\'emontr\'ees.

Au \S \ref{prem} nous donnons les cons\'equences les plus imm\'ediates du
th\'eor\`eme \ref{T1}, dont les corollaires \ref{T3} et \ref{T5bis}. Au
\S \ref{3}, nous introduisons la cohomologie de Lich\-ten\-baum (\`a
coefficients dans les complexes de cycles de Bloch),
\'etendons les classes de cycles motiviques
$l$-adiques de
\cite{glr} ($l\ne p=\car k$) en des classes \'emanant de la cohomologie
de Lichtenbaum motivique et cons\-trui\-sons une classe de cycle
$p$-adique correspondante; nous d\'emontrons ensuite leur bijectivit\'e
pour $X\in B_\tate(k)$. Nous en d\'eduisons le
corollaire
\ref{T2} et des propri\'et\'es de finitude pour la cohomologie motivique
de Lichtenbaum, meilleures que pour la cohomologie motivique \'etale (et
conjectur\'ees par Lichtenbaum dans \cite[introduction]{licht2} pour
toute vari\'et\'e projective lisse). Nous en d\'eduisons aussi une formule
pour le conoyau de l'application ``classe de cycle
$l$-adique enti\`ere"
$CH^n(X)\otimes\Z_l\to H^{2n}_\cont(X,\Z_l(n))$ (corollaire \ref{p5.7}).

Au \S \ref{ouv}, nous d\'emontrons quelques r\'esultats sur les ouverts
de $X$ et sur son corps des fonctions, dont
les corollaires \ref{T4} et \ref{T5}. Au \S \ref{retour}, nous tirons
la substantifique m\oe lle des arguments pr\'ec\'edents dans le cadre
motivique. Pour ce faire, nous introduisons une
cat\'egorie de ``motifs de Chow \'etales" $\Chow_\et(k)$, fort pratique
(voir
\cite{mil-ram} pour une approche diff\'erente, en termes de
r\'ealisations). Nous en d\'eduisons un curieux ``principe d'identit\'e
de Manin arithm\'etique" (corollaire \ref{idmanar}). Nous y d\'emontrons
aussi le corollaire \ref{T6}. Dans le
\S \ref{pros}, nous investigons la mesure dans laquelle les m\'ethodes
de cet article pourraient rapprocher l'\'ech\'eance de la
d\'e\-mons\-tra\-tion des conjectures indiqu\'ees. Enfin, dans
l'appendice, nous construisons une fonctorialit\'e covariante sur la
cohomologie motivique \'etale, n\'ecessaire pour la construction de
$\Chow_\et(k)$.

Une partie substantielle du travail de Kimura sur lequel nous nous
appuyons a \'et\'e obtenue ind\'ependamment par Peter J. O'Sullivan
(communication personnelle \`a Andr\'e et \`a l'auteur): notamment la
notion de dimension finie au sens de Kimura (semi-positivit\'e dans sa
terminologie) et le fait que le motif de Chow d'une vari\'et\'e
ab\'elienne est de dimension finie au sens de Kimura. Il n'obtient par
contre pas le th\'eor\`eme de nilpotence de Kimura, ni a fortiori la
proposition \ref{p0}.

Les lignes ci-dessus soulignent la dette que j'ai envers les travaux
an\-t\'e\-rieurs apparaissant dans la bibliographie: ils sont trop
nombreux pour \^etre cit\'es dans le d\'etail. En particulier, il aurait
sans doute \'et\'e possible de r\'ediger les \S\S \ref{3}--\ref{pros} en
termes de cohomologie motivique \'etale sans le travail de Geisser
\cite{geisser2}, mais les r\'esultats auraient
\'et\'e plus d\'esagr\'eables \`a \'enoncer et les d\'emonstrations plus
compliqu\'ees. 

Je remercie  Yves Andr\'e, Thomas Geisser, H\'el\`ene Esnault, Ofer
Gabber, Steven Lichtenbaum, Christophe Soul\'e et Eckart Viehweg pour
des commentaires pertinents sur la pr\'eparation de ce texte; en
particulier T. Geisser pour avoir attir\'e mon attention sur le
probl\`eme soulev\'e dans
\ref{meg1} et m'avoir incit\'e \`a \'ecrire l'appendice \ref{A}, et C.
Soul\'e pour m'avoir encourag\'e \`a ne pas me limiter aux produits de
courbes elliptiques. D'autre part, c'est avec Y. Andr\'e que j'ai
d\'egag\'e dans
\cite{ak} la proposition
\ref{p0} ci-dessous, r\'esultat cl\'e sur lequel repose tout ce travail.

\section{\'Equivalence rationnelle et \'equivalence
num\'erique}\label{ratnum}

Soient $\sA$ la cat\'egorie des motifs de Chow sur $k$ \`a coefficients
rationnels et $\bar \sA$ la cat\'egorie des motifs purs sur $k$ modulo
l'\'e\-qui\-va\-len\-ce num\'erique, \'egalement \`a coefficients
rationnels: cette derni\`ere est ab\'elienne semi-simple d'apr\`es
\cite{jannsen}. On a un foncteur plein 
\begin{align*}
\sA&\to \bar \sA\\
M&\mapsto \bar M.
\end{align*}

Pour tout objet $M$ de $\sA$ ou $\bar\sA$, on note $F_M$ l'endomorphisme
de Frobenius de $M$. Pour toute $k$-vari\'et\'e projective lisse $X$, on
note $h(X)$ (\resp $\bar h(X)$) le motif de $X$ dans $\sA$ (\resp son
image dans $\bar \sA$).

\begin{rem} Comme Voevodsky, nous adoptons la convention que le foncteur
$X\mapsto h(X)$ est \emph{covariant}, et non contravariant comme il est
d'u\-sa\-ge plus traditionnellement. Nous notons aussi $M\otimes L=M(1)$
au lieu de $M\otimes L=M(-1)$, o\`u $L$ est le motif de Lefschetz
($h(\P^1)=\un\oplus L$). Avec ces conventions, on a
$CH^n(X)\otimes\Q=\sA(h(X),L^n)$ et $CH_n(X)\otimes\Q=\sA(L^n,h(X))$ pour
toute $k$-vari\'et\'e projective lisse $X$, ainsi que $h(X)^\vee =
h(X)(-\dim X)$.
\end{rem}

Rappelons la d\'efinition suivante \cite{ki}:

\begin{defn} Un objet $M\in\sA$ est \emph{de dimension finie au sens de
Kimura} s'il existe une d\'ecomposition $M\simeq M_+\oplus M_-$ et deux
entiers $n_+,n_-\ge 0$ tels que
$\Lambda^{n_++1} M_+=\mathbf{S}^{n_-+1} M_-=0$.
\end{defn}

On a:

\begin{thm}[Kimura \protect{\cite[th. 4.2, cor. 5.11, prop.
6.9]{ki}}]\label{kim} La sous-ca\-t\'e\-go\-rie pleine de
$\sA$ form\'ee des motifs de dimenson finie au sens de Kimura est
additive, \'epaisse et rigide (stable par somme directe, facteurs
directs, produit tensoriel et dual). Elle contient les motifs des
vari\'et\'es ab\'eliennes.
\end{thm}

La proposition suivante est directement inspir\'ee de \cite[prop.
7.5]{ki}.

\begin{prop}[\protect{\cite[prop. 9.1.14]{ak}}]\label{p0} Soit $M\in\sA$
un motif de dimension finie au sens de Kimura. Alors le noyau de
$\sA(M,M)\to\bar\sA(\bar M,\bar M)$ est un id\'eal nilpotent.\qed
\end{prop}

\begin{lemme}\label{l1} Soit $N\in\sA$, de dimension finie au sens de
Kimura. Supposons que $\bar N$ soit simple et que $F_{\bar N}\ne 1$. Alors
$\sA(\un,N)=0$.
\end{lemme}

\begin{sloppypar}
\prf Soit $P\in\Q[T]$ le polyn\^ome minimal de
$F_{\bar N}$: comme $\bar\sA(\bar N,\bar N)$ est une $\Q$-alg\`ebre
\`a division, $P$ est un polyn\^ome irr\'eductible, dif\-f\'e\-rent de
$T-1$ par hypoth\`ese. D'apr\`es la proposition \ref{p0}, il existe
$n>0$ tel que $P(F_N)^n=0$. Soit $f\in \sA(\un,N)$. Alors $F_N f=f$,
d'o\`u $P(1)^nf=0$ et $f=0$.\qed
\end{sloppypar}

\begin{prop}\label{p1} Soient $M,M'\in\sA$ deux motifs de dimension finie
au sens de Kimura. Supposons que tout facteur simple $\bar N$ de
$(\bar M')^\vee\otimes \bar M$ v\'erifie soit $\bar N\simeq\un$, soit
$F_{\bar N}\ne 1$. Alors l'homomorphisme
$\sA(M',M)\to\bar\sA(\bar M',\bar M)$ est un isomorphisme.
\end{prop}

\prf Gr\^ace au th\'eor\`eme \ref{kim}, on se
ram\`ene par dualit\'e au cas o\`u $M'=\un$. En r\'eappliquant
la proposition \ref{p0}, on voit que toute d\'ecomposition de $1_{\bar
M}$ en somme d'idempotents orthogonaux se rel\`eve en une telle somme
dans $\sA(M,M)$. Pour cette d\'ecomposition, on a donc
\[M=\bigoplus M_i\]
o\`u  les $\bar M_i$ sont simples. De
plus, en utilisant toujours la proposition \ref{p0}, on voit que $\bar
M_i\simeq
\un\If M_i\simeq \un$.  La proposition \ref{p1} r\'esulte donc de
l'hypoth\`ese et du lemme \ref{l1}.\qed

\begin{rem} Le principe de la d\'emonstration du lemme \ref{l1} et de la
proposition \ref{p1} est tr\`es proche dans l'esprit de \cite[prop.
3.2]{geisser}.
\end{rem}

\begin{thm}\label{t1} Soient $X,X'$ deux $k$-vari\'et\'es projectives
lisses telles que $X\times_k X'\in B_\tate(k)$. Alors, pour tout
$n\in\Z$, l'homomorphisme
\begin{multline*}
CH^{\dim X+n}(X'\times X)\otimes\Q=\sA(h(X'),h(X)(n))\\
\to\bar\sA(\bar h(X'),\bar h(X)(n))=A^{\dim X+n}_\num(X'\times
X)\otimes\Q
\end{multline*} 
est un isomorphisme.
\end{thm}

\prf Comme ci-dessus, on se ram\`ene \`a $X'=\Spec k$. Remarquons que
$h(X)(n)$ est de dimension finie au sens de Kimura \cite[ex. 9.1]{ki}.
Par ailleurs, l'action de Frobenius sur les
$H^i(\bar X,\Q_l)$ est semi-simple: en effet, ils sont engendr\'es
multiplicativement par $H^1$ et cela r\'esulte alors des r\'esultats de
Weil \cite{weil}. En utilisant la conjecture de Tate, on en d\'eduit que
l'\'equivalence homologique et l'\'equivalence num\'erique
co\"{\i}ncident sur $X$ 
\cite[th. 2.9]{Tate}.  Il en r\'esulte alors (\cf
\cite[prop. 2.6]{milne}) que pour tout facteur simple $\bar N$ de $\bar
h(X)(n)$, on a soit $\bar N\simeq
\un$ soit $F_{\bar N}\ne 1$. L'hypoth\`ese de la proposition \ref{p1} est
donc v\'erifi\'ee, d'o\`u le th\'eor\`eme \ref{t1}.\hfill $\square$

\section{Premi\`eres applications}\label{prem}

\begin{cor}\label{c1} Soit $X\in B_\tate(k)$. Alors, pour tout $n$,
l'ordre du p\^ole de la fonction z\^eta
$\zeta(X,s)$ en $s=n$ est \'egal au rang de $CH^n(X)$ (qui est fini
d'apr\`es le th\'eor\`eme \ref{t1}).
\end{cor}

\prf Un cas particulier du th\'eor\`eme \ref{t1} est que, sur $X$,
l'\'e\-qui\-va\-len\-ce homologique (relative disons \`a la cohomologie
$l$-adique pour un nombre premier $l$ diff\'erent de la caract\'eristique
de $k$) est \'egale \`a l'\'equivalence
num\'erique. L'\'enonc\'e r\'esulte donc du th\'eor\`eme de Spie\ss, du
th\'eor\`eme \ref{t1} et de \cite[th. 2.9]{Tate}.\qed

La seconde application r\'esulte de Geisser \cite[th. 3.3]{geisser}:

\begin{cor}\label{t2} Si $X\in B_\tate(k)$, alors la conjecture de
Beilinson-Parshin est vraie pour $X$:
$K_i(X)$ est un groupe de torsion pour tout $i>0$.\qed
\end{cor}

Notons $H^i(X,\Q(n))$ les groupes de cohomologie motivique de $X$ \`a
coefficients rationnels: dans cette section, ils peuvent \^etre d\'efinis
\`a la Beilinson comme espaces propres de groupes de $K$-th\'eorie pour
les op\'erations d'Adams \cite{soule2}.

\begin{cor}\label{c2} Soient $X\in B_\tate(k)$, $d=\dim X$, $K=k(X)$ et
$n$ un entier $\ge 0$. Alors la conjecture de Beilinson-Soul\'e vaut pour
$K$ en poids $n$, et on a m\^eme
\[H^i(K,\Q(n))=0\text{ pour } i< n\]
dans les deux cas suivants:
\begin{thlist}
\item $d\le 2$
\item $n\le 2$.
\end{thlist}
De plus, sous la condition (i), $K_n^M(K)$ et $K_n(K)$ sont de
torsion pre\-mi\`e\-re \`a $p$ pour
$n\ge 3$ et $K_n^M(K)$ est de torsion impaire pour $n\ge 4$.
\end{cor}

\begin{sloppypar}
\prf \cf \cite[th. 3.4]{geisser}: on utilise le corollaire \ref{t2} et la
suite spectrale de coniveau pour la cohomologie motivique. Si $d\le 2$,
les corps r\'esiduels $F$ de $X$ autres que $K$ sont de dimension $0$ ou
$1$, donc les \'enonc\'es ``$H^i(F,\Q(n))=0$ pour $i< n$" et
``$H^i(F,\Q(n))=0$ pour $n>\dim F$" sont connus par les th\'eor\`emes de
Quillen \cite{quillen,quillen2} et de Harder \cite{harder}. Si $n\le 2$,
les poids intervenant dans les termes
$E_1^{p,q}$ avec $p>0$ sont $\le 1$, et l'\'enonc\'e est encore connu
\cite[cor. 6.8]{kratzer}. L'assertion ``premi\`ere \`a $p$" r\'esulte de
\cite{geilevp}, et la derni\`ere assertion pour $\car k\ne 2$ r\'esulte
de la conjecture de Milnor \cite{voem}.\hfill $\square$
\end{sloppypar}

\begin{cor}\label{c4.1} Soit $A$ un anneau de valuation discr\`ete de
corps des fractions $E$ et de corps r\'esiduel $K$, o\`u $K=k(X)$ avec
$X\in B_\tate(k)$, $\dim X\le 2$. Alors, pour tout
$n\ge 0$, la suite
\[0\to K_{n+1}(A)\to K_{n+1}(E)\to K_n(K)\to 0\]
est exacte.
\end{cor}

\prf Soit $m>0$
premier \`a $p=\car k$. Le diagramme commutatif au signe pr\`es
\[\begin{CD}
K_{n+2}(E,\Z/m)&\Surj&{}_mK_{n+1}(E)\\
@VVV @VVV\\
K_{n+1}(K,\Z/m)&\Surj&{}_m K_n(K)
\end{CD}\]
et la conjecture de Gersten pour $K_*(-,\Z/m)$ (Gillet
\cite{gillet}) montrent que l'homomorphisme r\'esidu
${}_mK_{n+1}(E)\to {}_mK_n(K)$. L'assertion r\'esulte donc du corollaire
\ref{c2} pour $n\ge 3$, et est bien connue pour $n\le 2$ puisqu'alors
$K_n(K)=K_n^M(K)$.\qed

\section{Classes de cycle motiviques}\label{3}

\begin{sloppypar}
\subsection{Cohomologie motivique} Nous d\'efinirons la cohomologie
motivique d'une $k$-vari\'et\'e lisse $X$ comme
l'hypercohomologie de Za\-ris\-ki $H^i_\Zar(X,\Z(n))$ des complexes de
cycles de Bloch d\'ecal\'es (\cf \cite{bloch,geilev,glr} pour les
d\'etails). Nous aurons aussi besoin de consid\'erer l'hypercohomologie
de ces complexes par rapport \`a d'autres topologies, et de les comparer.
En par\-ti\-cu\-lier, le r\'esultat ci-dessous sera d'usage constant:
\end{sloppypar}

\begin{prop}[\protect{\cite[prop. 1.18]{glr}}]\label{p4} Pour toute
vari\'et\'e lisse
$X$, les homomorphismes
\[H^i_\Zar(X,\Q(n))\to H^i_\et(X,\Q(n))\]
sont des isomorphismes.\qed
\end{prop}

On a aussi:

\begin{prop}\label{p4bis} Soit $X\in B_\tate(k)$. Alors
$H^i_\Zar(X,\Q(n))=0$ pour $i\ne 2n$.
\end{prop}

\prf Notons que l'\'enonc\'e est trivialement vrai pour $i>2n$ et toute 
vari\'et\'e lisse $X$ (sur un corps quelconque) par d\'efinition de la
cohomologie motivique. Pour $i<2n$, cela r\'esulte du corollaire \ref{t2}
et du fait que $H^i_\Zar(X,\Q(n))\simeq gr_\gamma^n K_{2n-i}(X)$
\cite{bloch}. (Variante: raisonner directement comme dans la
d\'emonstration du th\'eor\`eme \ref{t1}, en utilisant le fait que
Frobenius agit sur ce groupe par multiplication par $q^n$, \cf
\cite{geisser}.)\hfill $\square$

\subsection{Cohomologie de Lichtenbaum}\label{cL}

Nous allons utiliser la cohomologie introduite par Lichtenbaum dans
\cite{licht2}\footnote{Lichtenbaum d\'enomme cette cohomologie
\emph{Weil-\'etale cohomology}. Il nous para\^{\i}t plus pratique et plus
juste de la rebaptiser cohomologie de Lichtenbaum, la terminologie
``cohomologie de Weil" pouvant de toute fa\c con pr\^eter \`a
confusion.}. Cette cohomologie a \'et\'e \'egalement \'etudi\'ee par
Geisser \cite{geisser2}. 

Rappelons \cite[prop. 2.2]{licht} que la cat\'egorie des faisceaux pour la
topologie de Lichtenbaum sur une $k$-vari\'et\'e $X$ de type fini est
\'equivalente \`a la cat\'egorie des faisceaux
\'etales $\Z$-\'equivariants sur $\bar X=X\times_k \bar k$, o\`u l'action
de $\Z$ sur $\bar X$ est donn\'ee par le morphisme de Frobenius
galoisien. Pour un tel faisceau $\sF$, on pose
\[H^0_W(X,\sF)=H^0_\et(\bar X,\sF)^\Z\]
et la cohomologie de Lichtenbaum $H^*_W$ est d\'efinie comme les foncteurs
d\'eriv\'es de $\sF\mapsto H^0_W(X,\sF)$. En particulier on a une suite
spectrale \cite[prop. 2.3]{licht2}
\begin{equation}\label{e16}
H^p(\Z,H^q_\et(\bar X,C))\Rightarrow
H^{p+q}_W(X,C)
\end{equation}
pour tout complexe de faisceaux \'etales $\Z$-\'equivariants.

Pour comparer cohomologie \'etale et cohomologie de Lichtenbaum, le point
de vue de Geisser \cite[\S 2]{geisser2} est d'identifier la cat\'egorie
des faisceaux \'etales sur $X$ \`a la cat\'egorie des faisceaux \'etales
$\hat{\Z}$-\'equivariants (topologiques discrets) sur $\bar X$, o\`u
$\hat{\Z}=Gal(\bar k/k)$
\cite[Exp. XIII, 1.1.3]{deligne3}. Soit $\sT_\Z(X)$ (\resp
$\sT_{\hat{\Z}}(X)$) la cat\'egorie des $\bar X$-faisceaux \'etales
$\Z$-\'equivariants (\resp $\hat{\Z}$-\'e\-qui\-va\-riants): on a des
foncteurs adjoints \'evidents
\[\gamma_X^*:\sT_{\hat\Z}(X)\to
\sT_{\Z}(X),\quad
(\gamma_X)_*:\sT_\Z(X)\to \sT_{\hat{\Z}}(X).\]

\begin{sloppypar}
Notons $\sT_\Z$ et $\sT_{\hat\Z}$ les cat\'egories de faisceaux
correspondantes sur les ``grands sites lisses" correspondant \`a la
cat\'egorie $Sm/k$ des $k$-vari\'et\'es lisses de type fini. On a des
foncteurs adjoints correspondants
\[\gamma^*:\sT_{\hat\Z}\to\sT_{\Z},\quad
\gamma_*:\sT_\Z\to \sT_{\hat{\Z}}\]
induisant des foncteurs adjoints sur les cat\'egories d\'eriv\'ees
\[\gamma^*:D(\sT_{\hat\Z})\to D(\sT_{\Z}),\quad
R\gamma_*:D(\sT_\Z)\to D(\sT_{\hat{\Z}}).\]
\end{sloppypar}

\begin{sloppypar}
\subsubsection{Mise en garde} \label{meg1} Le foncteur $R\gamma_*$ n'est
pas conservatif. Par exemple, soit $\Q\langle n\rangle$ le $\Z[\Z]$-module
de support
$\Q$, l'action du g\'en\'erateur de $\Z$ \'etant donn\'ee par $r\mapsto
q^n r$. Alors
$\Q\langle n\rangle$ d\'efinit un objet de $\sT_\Z$, mais  on a
$H^*_W(X,\Q\langle n\rangle)=0$ pour tout $X\in Sm/k$ d\`es que $n\ne
0$: c'est imm\'ediat \`a partir de \eqref{e16}.

Pour obvier cet inconv\'enient, introduisons la sous-cat\'egorie \'epaisse
$D_{in}(\sT_\Z)$ de $D(\sT_\Z)$ form\'ee des complexes $C$ tels que
$H^*_W(X,C)=0$ pour tout $X\in Sm/k$, et notons
\[\bar D(\sT_\Z)=D(\sT_\Z)/D_{in}(\sT_\Z).\]
\end{sloppypar}

Soit $\bar\gamma^*$ le compos\'e de $\gamma^*$ avec le foncteur de
localisation $D(\sT_\Z)\to \bar D(\sT_\Z)$. Alors $R\gamma_*$ se
factorise en un foncteur exact
\[\bar R\gamma_*:\bar D(\sT_\Z)\to D(\sT_{\hat\Z})\]
qui est conservatif par construction, et adjoint \`a droite de
$\bar\gamma^*$. Nous aurons besoin de 
$\bar R\gamma_*$ plut\^ot que de $R\gamma_*$ pour une formulation
correcte du th\'eor\`eme \ref{t6} (iv).

\begin{sloppypar}
\subsection{Cohomologie de Lichtenbaum motivique} Nous
con\-si\-d\'e\-re\-rons prin\-ci\-pa\-le\-ment la cohomologie de
Lichtenbaum
\`a coefficients dans les complexes de cycles de Bloch. Pour tout $X\in
Sm/k$, le foncteur
$R(\gamma_X)_*$ induit des homomorphismes canoniques
\[H^i_\et(X,\Z(n))\to H^i_W(X,\Z(n)).\]
\end{sloppypar}

Le r\'esultat suivant de Geisser sera d'usage constant: notons $\beta$ le
Bockstein et $e$ le g\'en\'erateur de $H^1_W(k,\Z)=Hom(\Z,\Z)$ envoyant
(par exemple) le Frobenius g\'eom\'etrique sur $1$. On a de longues
suites exactes
\cite[th. 6.1]{geisser2}
\begin{multline}\label{e1}
\dots\to H^i_\et(X,\Z(n))\to H^i_W(X,\Z(n))\\
\to H^{i-1}_\et(X,\Z(n))\otimes\Q\stackrel{\partial}{\to}
H^{i+1}_\et(X,\Z(n))\to\dots
\end{multline}
o\`u $\partial$ est donn\'ee par la composition
\begin{multline*}
H^{i-1}_\et(X,\Z(n))\otimes \Q=H^{i-1}_\et(X,\Q(n))\to
H^{i-1}_\et(X,\Q/\Z(n))\overset{\cdot e}{\to}\\ 
H^{i}_\et(X,\Q/\Z(n))\overset{\beta}{\to}H^{i+1}_\et(X,\Z(n))
\end{multline*}
(\cf \cite[prop. 9.12]{tate}).

On a aussi:

\begin{prop}[\protect{\cite[th. 6.5 et 6.7]{geisser2}}]\label{p2} a)
Pour toute
$k$-vari\'et\'e lisse $X$ de dimension
$d$ et pour tout $n\ge 0$, on a $H^i_W(X,\Z(n))=0$ pour
$i>\sup(2d+1,n+d+1)$.\\ 
b) Si $X$ est de plus projective, on a un
diagramme commutatif 
\[\begin{CD}H^{2d}_W(X,\Z(d))@>>> H^{2d}_\et(\bar X,\Z(d))^\Z@= \Z\oplus
(fini)\\
@V{\cdot e}VV @VVV @VVV\\
H^{2d+1}_W(X,\Z(d))@<\sim<< H^{2d}_\et(\bar X,\Z(d))_\Z@= \Z
\end{CD}\]
o\`u l'isomorphisme $H^{2d}_\et(\bar X,\Z(d))_\Z= \Z$ provient de
l'homomorphisme ``de\-gr\'e"
\[H^{2d}_\et(\bar X,\Z(d))\simeq CH^d(\bar X)\to \Z\]
(\cf \cite[d\'em. du lemme 6.4 a)]{geisser2}).\qed
\end{prop}

Nous aurons aussi besoin de la

\begin{prop}\label{p3} Pour toute vari\'et\'e lisse $X$, $H^i_W(X,\Z(n))$
est de torsion pour $i>2n+1$.
\end{prop}

\begin{sloppypar}
\prf Cela r\'esulte par exemple du fait que $H^i_\Zar(\bar
X,\Q(n))= H^i_\et(\bar X,\Q(n))=0$ pour $i>2n$ et de la suite spectrale
\eqref{e16}.\hfill $\square$
\end{sloppypar}

\subsection{Classe de cycle $l$-adique}\label{l} Fixons un nombre premier
$l\ne p$. Nous allons utiliser une application \emph{classe de cycle
motivique
$l$-adique pour la cohomologie de Lichtenbaum}. Cela revient \`a
\'etendre les homomorphismes de \cite[\S 1.4]{glr}
\[H^i_\et(X,\Z(n))\otimes\Z_l\to H^i_\cont(X,\Z_l(n))\]
en des homomorphismes provenant de $H^i_W(X,\Z(n))\otimes\Z_l$.

Indiquons rapidement comment on proc\`ede: notons $\alpha:(Sm/k)_\et\to
(Sm/k)_\Zar$ le foncteur de
projection entre les sites d\'efinis par la cat\'egorie $Sm/k$ des
$k$-vari\'et\'es lisses de type fini munie respectivement des topologies
\'etale et de Zariski.  Dans \cite[th. 1.17]{glr}, nous
avons d\'efini un objet $\Z(n)\in D^{-}(Ab((Sm/k)_\Zar))$ dont la
res\-tric\-tion
\`a toute $k$-vari\'et\'e lisse $X$ est quasi-isomorphe au complexe de
cycles de Bloch de poids $n$ sur $X$, puis un morphisme dans la
cat\'egorie d\'eriv\'ee \cite[(1.8)]{glr}
\begin{equation}\label{e5}
\alpha^*\Z(n)\oo^L\Z_l\to\Z_l(n)^c_\et
\end{equation}
o\`u $\Z_l(n)^c_\et=R\plim\mu_{l^\nu}^{\otimes n}$ dans la
cat\'egorie d\'eriv\'ee des faisceaux \'etales, que nous identifions
canoniquement \`a
$D(\sT_{\hat \Z})$ (\cf \ref{cL}). 

D\'efinissons de m\^eme
$\Z_l(n)^c_W=R\plim\gamma^*\mu_{l^\nu}^{\otimes n}$ dans $D(\sT_\Z)$
(ibid.). On a alors un morphisme compos\'e
\begin{equation}\label{e3}
\gamma^*\alpha^*\Z(n)\to \gamma^*\Z_l(n)^c_\et\to \Z_l(n)^c_W.
\end{equation}

C'est le morphisme cherch\'e: il induit des homomorphismes
\[H^i_W(X,\Z(n))\otimes\Z_l\to H^i_W(X,\Z_l(n)^c_W).\]

Il reste \`a remarquer que l'adjoint $\Z_l(n)^c_\et\to R\gamma_*
\Z_l(n)^c_W$ du morphisme $\gamma^*\Z_l(n)^c_\et\to \Z_l(n)^c_W$ est un
quasi-isomorphisme: cela r\'esulte du fait que $R\gamma_*$ et
$R\plim$ commutent et de l'analogue de \eqref{e1} \`a coefficients
$\mu_{l^\nu}^{\otimes n}$, qui montre que $\mu_{l^\nu}^{\otimes n}\to
R\gamma_*\gamma^*\mu_{l^\nu}^{\otimes n}$ est un quasi-isomorphisme pour
tout $n$ (\cf \cite[cor. 3.4]{geisser2}). On en conclut que
l'homomorphisme canonique
\[H^i_\cont(X,\Z_l(n))\to H^i_W(X,\Z_l(n)^c_W)\]
est bijectif, et on en d\'eduit bien l'application ``classe de cycle"
promise:
\begin{equation}\label{e4}
H^i_W(X,\Z(n))\otimes\Z_l\to H^i_\cont(X,\Z_l(n)).
\end{equation}

(Bien entendu, on pourrait aussi proc\'eder na\"{\i}vement \`a partir de
\cite[\S 3.7]{geilev}, voir aussi \cite[\S 4]{bloch3}, qui est de toute
fa\c con \`a la base de notre
construction, en d\'efinissant \eqref{e4} ``vari\'et\'e par vari\'et\'e"
et en \'evitant \cite[(1.8)]{glr}, mais les d\'etails seraient plus
p\'enibles \`a r\'ediger.)

En fait:

\begin{lemme}\label{l2} L'homomorphisme \eqref{e4} n'est autre que celui
induit par le morphisme $\Phi_n:\alpha^*\Z(n)\oo\limits^L\Z_l(0)_\et^c\to
\Z_l(n)^c_\et$ de
\cite[conj. 3.2]{glr}.
\end{lemme}

\prf Rappelons que $\Phi_n$ est d\'efini comme la composition
\[\alpha^*\Z(n)\oo^L\Z_l(0)^c_\et\to\Z_l(n)^c_\et\oo^L\Z_l(0)^c_\et\to
\Z_l(n)^c_\et\] 
o\`u la premi\`ere fl\`eche se d\'eduit de \eqref{e5} et la suivante est
donn\'ee par le cup-produit. Nous allons voir que $\Phi_n$ se d\'eduit de
\eqref{e3} par application du foncteur $R\gamma_*$.

En effet, on a d\'ej\`a vu le quasi-isomorphisme $\Z_l(n)^c_\et\iso
R\gamma_*\Z_l(n)^c_W$. D'autre part, d'apr\`es \cite[th. 3.3]{geisser2},
le morphisme de double adjonction
\begin{equation}\label{e14}
R\gamma_*\Z\oo^L\alpha^*\Z(n)\to
R\gamma_*\gamma^*\alpha^*\Z(n)
\end{equation}
est \'egalement un quasi-isomorphisme.
Dans \cite[th. 4.6 et 6.3]{tate}, nous avons \'etabli un
quasi-isomorphisme
\[\Z_l(0)^c_\et\simeq \Z^c\otimes\Z_l\]
o\`u $\Z^c$ est un certain complexe de longueur $1$. Enfin, dans
\cite[th. 4.2]{geisser2}, Geisser \'etablit un quasi-isomorphisme
$\Z^c\simeq R\gamma_*\Z$. L'assertion r\'esulte ai\-s\'e\-ment de tous ces
quasi-isomorphismes, en suivant leurs d\'efinitions respectives.\hfill
$\square$

\subsection{Classe de cycle $p$-adique}\label{p} Posons
\[\Z_p(n)^c_\et=R\plim \nu_{r}(n)[-n]\quad (n\ge 0)\]
o\`u $\nu_{r}(n)$ est le $n$-i\`eme faisceau de Hodge-Witt logarithmique
de niveau $r$, et de m\^eme
\[\Z_p(n)^c_W=R\plim \gamma^*\nu_{r}(n)[-n].\]

Pour $n<0$, on pose $\Z_p(n)^c_\et=0$ et $\Z_p(n)^c_W=0$.

En passant \`a la limite sur la construction de \cite[d\'em. du
th. 8.3]{geilevp}, on obtient une classe de cycle motivique $p$-adique
pour la cohomologie \'etale
\[\alpha^*\Z(n)\otimes\Z_p\to \Z_p(n)^c_\et\]
puis une classe de cycle motivique $p$-adique pour la cohomologie de
Lichtenbaum
\begin{equation}
\gamma^*\alpha^*\Z(n)\otimes\Z_p\to \gamma^*\Z_p(n)^c_\et\to \Z_p(n)^c_W
\end{equation}
et plus concr\`etement, pour une vari\'et\'e lisse $X$ donn\'ee, des
homomorphismes
\begin{equation}\label{e4p}
H^i_W(X,\Z(n))\otimes\Z_p\to H^i_\cont(X,\Z_p(n))
\end{equation}
o\`u $H^i_\cont(X,\Z_p(n)):=H^i_\et(X,\Z_p(n)^c_\et)\iso
H^i_W(X,\Z_p(n)^c_W)$ (\cf \ref{l}).

\subsection{Bijectivit\'e des classes de cycle}

\begin{thm}\label{t4} Si $X\in B_\tate(k)$,
\eqref{e4} est un isomorphisme pour tous $l\ne p$, $i,n$, et \eqref{e4p}
est un isomorphisme pour tous $i,n$.
\end{thm}

\prf Traitons d'abord le cas de \eqref{e4}. D'apr\`es le lemme \ref{l2},
il suffit de voir que le morphisme induit par $\Phi_n$ est un isomorphisme
pour tout $n$; cela r\'esulte du th\'eor\`eme \ref{t1}, du corollaire
\ref{c1} et de \cite[th. 3.4]{glr}. (On pourrait aussi proc\'eder
directement \`a partir de la conjecture de Tate et du th\'eor\`eme
\ref{t1}, mais ce serait plus laborieux.) On en d\'eduit d\'ej\`a:

\begin{lemme}\label{l7} Pour $X\in B_\tate(k)$, on a $H^i_W(X,\Q(n))=0$
pour $i<2n$ et la composition
\[CH^n(X)\otimes\Q\simeq H^{2n}_\Zar(X,\Q(n))\to H^{2n}_W(X,\Q(n))\]
est un isomorphisme pour tout $n\ge 0$.
\end{lemme}

\begin{sloppypar}
\prf Le premier point
r\'esulte par r\'ecurrence sur $i$ de la suite exacte \eqref{e1} (\cf
propositions \ref{p4} et \ref{p4bis}). Cette m\^eme suite exacte implique
alors que $H^{2n}_\et(X,\Q(n))\to H^{2n}_W(X,\Q(n))$ est un isomorphisme,
et on conclut  par la proposition \ref{p3}. \qed
\end{sloppypar}

D\'emontrons maintenant que \eqref{e4p} est aussi un isomorphisme. On
proc\`ede directement, en imitant la m\'ethode de \cite[\S 3.4]{glr}: on
remarque d'abord que la version de
\eqref{e4p} \`a coefficients $\Q_p/\Z_p$ est un isomorphisme d'apr\`es
\cite{geilevp} (ce fait est vrai pour toute $k$-vari\'et\'e lisse). Par le
lemme des 5, on se ram\`ene \`a d\'emontrer que \eqref{e4p}$\otimes\Q$ est
un isomorphisme. 

On sait d\'ej\`a que $H^i_W(X,\Q(n))=0$ pour $i<2n$ (par le lemme
\ref{l7}) et pour $i>2n+1$ (par la proposition \ref{p3}). La m\^eme chose
est vraie pour les groupes $H^i_\cont(X,\Q_p(n))$, en utilisant un
th\'eor\`eme de Gabber \cite[th. 3]{ctss} et la suite exacte longue
\begin{multline*}
\dots\to H^{i-n-1}_\et(X,\nu_{\infty}(n))\to H^i_\cont(X,\Z_p(n))\to\\
H^i_\cont(X,\Q_p(n))\to H^{i-n}_\et(X,\nu_{\infty}(n))\to\dots
\end{multline*}

\begin{sloppypar}
Il reste donc \`a traiter les cas $i=2n$ et $i=2n+1$. Milne
\cite[prop. 5.4]{milne2} a d\'emontr\'e qu'on peut utiliser l'action
galoisienne sur les groupes $H^i(\bar X,\Q_p(n))$ \`a la place de
l'action de Frobenius sur la cohomologie cris\-tal\-li\-ne pour calculer
la fonction
$\zeta(X,s)$. En particulier, en tenant compte de l'hypoth\`ese de Riemann
\cite{deligne} et de \cite{km}, cela d\'emontre que
\begin{equation}\label{e13}
H^i(\bar X,\Q_p(n))^G=H^i(\bar X,\Q_p(n))_G=0\text{ pour }
i\ne 2n.
\end{equation}
\end{sloppypar}

Consid\'erons le diagramme commutatif (\cf \cite[d\'em. de la prop.
3.9]{glr})
\begin{equation}
\begin{CD}\label{e10}
H^{2n}_W(X,\Z(n))\otimes\Q_p&\to&
H^{2n}_\cont(X,\Q_p(n))&\overset{\sim}{\To}&
H^{2n}_\cont(\bar X,\Q_p(n))^G\\
@V{\cdot e}V{\wr}V @V{\cdot e}VV @V{f}VV\\
H^{2n+1}_W(X,\Z(n))\otimes\Q_p&\to&
H^{2n+1}_\cont(X,\Q_p(n))&\overset{\sim}{\longleftarrow}&
H^{2n}_\cont(\bar X,\Q_p(n))_G
\end{CD}
\end{equation}
o\`u  $f$ est la composition
\[H^{2n}_\cont(\bar X,\Q_p(n))^G\Inj H^{2n}_\cont(\bar X,\Q_p(n))\Surj
H^{2n}_\cont(\bar X,\Q_p(n))_G\]
(voir \cite[prop. 6.5]{milne2} pour la commutativit\'e du carr\'e de
droite). La fl\`eche verticale de gauche est un isomorphisme par le cas
$l\ne p$, et les fl\`eches horizontales de droite sont des isomorphismes
gr\^ace
\`a
\eqref{e13}.

En proc\'edant comme dans \cite[\S 2]{Tate} et en utilisant \cite[prop.
5.4]{milne2} et \cite{km}, on voit que le corollaire \ref{c1} implique que
$f$ est bijective et que la composition
\[CH^n(X)\otimes\Q_p\to H^{2n}_W(X,\Z(n))\otimes\Q_p\to
H^{2n}_\cont(\bar X,\Q_p(n))^G\]
est surjective. De plus, le th\'eor\`eme \ref{t1} implique que cette
composition est injective. En utilisant le lemme \ref{l7}, on en d\'eduit
que la composition horizontale sup\'erieure dans \eqref{e10} est
bijective. Il en r\'esulte que tous les homomorphismes de \eqref{e10}
sont bijectifs.\hfill$\square$

\subsection{Cons\'equences}

\begin{cor}\label{c3} Si $X\in B_\tate(k)$ et $d=\dim X$,\\ 
a) L'accouplement
\begin{equation}\label{e2}
H^{2n}_W(X,\Z(n))\times H^{2d-2n}_W(X,\Z(d-n))\to
H^{2d}_W(X,\Z(d))\to \Z
\end{equation}
est parfait modulo torsion pour tout $n$.\\
b) Pour tout $(i,n)$, l'accouplement
\begin{multline*}
H^i_W(X,\Z(n))_\tors\times H^{2d+1-i}_W(X,\Q/\Z(d-n))\\
\to H^{2d+1}_W(X,\Q/\Z(d))\iso \Q/\Z
\end{multline*}
induit un accouplement parfait de groupes finis
\[H^i_W(X,\Z(n))_\tors\times H^{2d+2-i}_W(X,\Z(d-n))_\tors\to \Q/\Z.\]
c) (\cf \cite[introduction]{licht2}) Les groupes $H^i_W(X,\Z(n))$ sont de
type fini pour tout $i$, finis
pour $i\notin\{2n,2n+1\}$ et nuls pour $i\le 0$ (si $n>0$).\\
d) Le noyau et le conoyau du cup-produit par $e$
\[H^{2n}_W(X,\Z(n))\to H^{2n+1}_W(X,\Z(n))\]
sont finis.\\
e) L'homomorphisme canonique
\[H^i_\et(X,\Z(n))\to H^i_W(X,\Z(n))\]
est un isomorphisme pour $i\le 2n$.
\end{cor}

\prf Il r\'esulte de \cite{deligne,km,gabber} que
$H^i_\cont(X,\Z_l(n))_\tors$ est fini pour tout $i\in\Z$, y compris $l=p$,
et nul pour presque tout $l$ (\cf \cite[th. 2 et 3]{ctss}). 
De plus, $H^i_\cont(X,\Z_l(n))$ est fini pour $i\ne 2n,2n+1$ (ibid.). Il
r\'esulte d'abord de ceci et du th\'eor\`eme \ref{t4} que c) est vrai
pour $i\notin\{2n,2n+1\}$, et est vrai pour $i\in\{2n,2n+1\}$ apr\`es
tensorisation par $\Z_l$ (pour tout $l$).

L'assertion b) est vraie apr\`es tensorisation par $\Z_l$ par le cas
de la cohomologie \'etale continue (\cf \cite[th. 1.14]{milne2} pour le
cas
$l=p$); sa version enti\`ere en r\'esulte directement. De m\^eme, a) est
vrai apr\`es tensorisation par
$\Z_l$. La version enti\`ere de a) et de c) r\'esulte alors du lemme
suivant, que nous n'avons pas trouv\'e dans la litt\'erature:

\begin{lemme}\label{l3} Soient $R$ un anneau de Dedekind et $A\times B\to
R$ un accouplement de deux $R$-modules $A,B$ sans torsion.\\
a) Supposons que cet accouplement devienne parfait apr\`es tensorisation
par $R_l$ pour tous les id\'eaux maximaux $l$ de $R$, o\`u $R_l$
d\'esigne le compl\'et\'e de $R$ en $l$. Alors il est parfait.\\
b) Soit $K$ le corps des fractions de $R$. Si, de plus, $A\otimes K$
ou $B\otimes K$ est un $K$-espace vectoriel de dimension finie, alors $A$
et $B$ sont de type fini.
\end{lemme}

\prf a) Soit $R_{(l)}$ le localis\'e de $R$ en $l$ . Comme l'extension
$R_l/R_{(l)}$ est fid\`element plate, l'hypoth\`ese vaut en rempla\c cant
$R_l$ par $R_{(l)}$. Consid\'erons  l'homomorphisme  (injectif)
\begin{equation}\label{e6}
A\to Hom(B,R).
\end{equation}

On a une cha\^{\i}ne d'homomorphismes
\[A_{(l)}\to Hom(B,R)_{(l)}\to Hom(B_{(l)},R_{(l)}).\]

Il est facile de voir que le second homomorphisme est injectif. De plus,
la composition est bijective par hypoth\`ese. Par cons\'equent, le
premier homomorphisme est lui aussi bijectif. Le conoyau $M$ de
\eqref{e6} v\'erifie donc $M_{(l)} = 0$ pour tout $l$. Il en r\'esulte que
$M = 0$ et que l'accouplement est parfait.

b) Pour fixer les id\'ees, supposons que $\dim_K B\otimes K<\infty$. Soit
$(b_1,\dots, b_n)$ une famille  d'\'el\'ements de
$B$ qui d\'efinit une base de $B\otimes K$, et soit $B'$ le sous-module
de $B$ engendr\'e par les $b_i$. Comme pour tout $b\in B$, il existe 
$r\in R$ tel que $rb\in B'$, l'application $Hom(B,R)\to Hom(B',R)$ est
injective. Donc $A$ s'injecte dans un $R$-module de type fini et il est
donc lui-m\^eme de type fini. En \'echangeant $A$ et $B$, on obtient la
m\^eme conclusion pour $B$.\qed

Revenons \`a la d\'emonstration du corollaire \ref{c3}: il reste \`a
d\'emontrer d) et e). L'assertion d) r\'esulte de c) et du fait qu'elle
est vraie apr\`es tensorisation par $\Q$, comme il r\'esulte de la
d\'emonstration du th\'eor\`eme \ref{t4}. Enfin, e) d\'ecoule de c) et de
\eqref{e1}, par r\'ecurrence sur $i$.\qed

\begin{cor}\label{c3l} Les conjectures de \cite[\S 7]{licht}
sont vraies pour $X\in B_\tate(k)$.
\end{cor}

\prf Rappelons tout d'abord ces conjectures \cite[\S 7]{licht} en termes
modernes; \'etant donn\'e une $k$-vari\'et\'e projective lisse $X$:

\begin{sloppypar}
\begin{enumerate}
\item $H^i_\et(X,\Z(n))=0$ pour $i$ grand.
\item $H^{2n}_\et(X,\Z(n))$ est un groupe ab\'elien de type fini.
\item $H^{i}_\et(X,\Z(n))$ est fini pour $i\ne 2n,2n+2$, nul pour $i\le
0$ lorsque $n>0$\footnote{Cette conjecture de nullit\'e r\'esulte de
l'axiome (1) de \protect\cite[\S 3]{licht}; il est pratique de l'ins\'erer
ici.}.
\item $H^{2d+2}_\et(X,\Z(d))$ est canoniquement isomorphe \`a $\Q/\Z$,
o\`u $d=\dim X$.
\item L'accouplement 
\[H^{i}_\et(X,\Z(n))\times H^{2d+2-i}_\et(X,\Z(d-n))\to
H^{2d+2}_\et(X,\Z(d))\iso \Q/\Z\]
est ``parfait", au sens qu'il d\'efinit une dualit\'e parfaite de
groupes finis  pour
$i\ne 2n$ et une dualit\'e parfaite entre un groupe de type fini et un
groupe de cotype fini pour $i=2n$. En particulier,
$\rg H^{2d}_\et(X,\Z(d))=1$.
\item Les groupes $H^{2n}_\et(X,\Z(n))$ et $H^{2d-2n}_\et(X,\Z(d-n))$ ont
le m\^eme rang $m(n)$.
\item $m(n)$ est l'ordre du p\^ole de la fonction $Z(X,t)$ en
$t=q^{-n}$ ($\zeta(X,s)=Z(X,q^{-s})$).
\item $\lim\limits_{t\to q^{-n}} (1-q^nt)^{m(n)}Z(X,t)=\pm
q^{\chi(X,\sO_X,n)}\chi(X,\Z(n))$, avec
\begin{multline*}\chi(X,\Z(n))=\\
\prod_{i\ne 2n,2n+2}
|H^i_\et(X,\Z(n))|^{(-1)^i}\cdot
\frac{|H^{2n}_\et(X,\Z(n))_\tors||H^{2n+2}_\et(X,\Z(n))_\cotors|}{R_n(X)}
\end{multline*}
o\`u $R_n(X)$ est la valeur absolue du d\'eterminant de l'accouplement 
\begin{multline*}
H^{2n}_\et(X,\Z(n))/\tors\times
H^{2d-2n}_\et(X,\Z(d-n))/\tors\\
\to  H^{2d}_\et(X,\Z(d))/\tors\iso\Z
\end{multline*}
par rapport \`a des bases quelconques de $H^{2n}_\et(X,\Z(n))/\tors$ et de
\break $H^{2d-2n}_\et(X,\Z(d-n))/\tors$, et
\[\chi(X,\sO_X,n)=\sum_{\substack{0\le i\le n\\ 0\le j\le d}}
(-1)^{i+j}(n-i){h_{ij}},\quad h_{ij}=\dim H^j(X,\Omega^i).\]
\end{enumerate}
\end{sloppypar}

En fait, $X$ \'etant quelconque, les \'enonc\'es 1 et 4 ne sont plus des
conjectures, ni le fait que  $\rg H^{2d}_\et(X,\Z(d))=1$: vu \eqref{e1},
le point a) de la proposition \ref{p2} implique la conjecture 1 de
Lichtenbaum pour
$i>2d+2$ et le point  b) implique la conjecture 4 et le fait que
$H^{2d}_\et(X,\Z(d))$ est de type fini et de rang $1$. Cela peut
d'ailleurs s'obtenir directement.

Il reste \`a traiter les conjectures 2, 3, 5, 6, 7 et 8. La conjecture 2
r\'esulte du corollaire \ref{c3} c) et e), ainsi que la conjecture 3 pour
$i<2n$; le reste d\'ecoule du corollaire \ref{c3} c) et de la suite
exacte \eqref{e1}. De m\^eme, 5 d\'ecoule du corollaire \ref{c3} a) et
b), via \eqref{e1}, et 6 d\'ecoule du corollaire \ref{c3} a) et e).
Enfin, 7 d\'ecoule du corollaire \ref{c3} e) et du corollaire
\ref{c1}.

Pour 8, on pr\'ef\`ere d\'emontrer la variante de \cite[th.
8.1]{geisser2} (voir aussi \cite[th. 10.7]{mil-ram}), qui lui est
\'e\-qui\-va\-len\-te via \eqref{e1}:

\begin{sloppypar}
\begin{itemize}
\item[$8_W$.] $\lim\limits_{t\to q^{-n}} (1-q^nt)^{m(n)}Z(X,t)=\pm
q^{\chi(X,\sO_X,n)}\chi(X,\Z(n))$, avec
\[\chi(X,\Z(n))=\prod_{i}
|H^i_W(X,\Z(n))_\tors|^{(-1)^i}\cdot R_n(X)^{-1}\]
et $\chi(X,\sO_X,n)$ comme dans 8, o\`u $R_n(X)$ est la valeur absolue du
d\'eterminant de l'accouplement
\eqref{e2} vu modulo torsion par rapport \`a des bases quelconques de
$H^{2n}_W(X,\Z(n))/\tors$ et de
$H^{2d-2n}_W(X,\Z(d-n))/\tors$.
\end{itemize}
\end{sloppypar}

On peut proc\'eder comme dans \cite[th. 4.3 et cor. 5.5]{milne3} (\cf
\cite[cor. 7.10 et th. 9.20]{tate} et \cite[d\'em. du th.
8.1]{geisser2}).\hfill $\square$

\begin{rem}\label{r1} Il est facile de pr\'eciser le signe dans la
conjecture 8 de Lichtenbaum: on remarque simplement que $\zeta(X,s)$ est
\`a valeurs r\'eelles positives pour $s$ r\'eel assez grand et ne s'annule
pas pour $s>d$. Le signe en $s=n$ ne d\'epend alors que de l'ordre des
z\'eros pour $s>n$: on trouve qu'il est \'egal \`a $(-1)^{\sum_{a>n}
m(a)}$.
\end{rem}
\enlargethispage*{20pt}

\begin{cor}\label{p5.7} Soit $X\in B_\tate(k)$. Alors, pour tout $l$, on
a une suite exacte
\begin{multline*}
0\to CH^n(X)\{l\}\to
CH^n(X)\otimes\Z_l\oplus
H^{2n-1}_\et(X,\Q_l/\Z_l(n))\\
\to H^{2n}_\cont(X,\Z_l(n))
\to
CH^n(X)\otimes\Q_l/\Z_l\to H^{2n}_\et(X,\Q_l/\Z_l(n)).\end{multline*}
En particulier, on a une suite exacte
\begin{multline*}
0\to \Ker(CH^n(X)\otimes\Q_l/\Z_l\overset{\bar{cl}}{\To}
H^{2n}_\et(X,\Q_l/\Z_l(n)))\\
\to
\Coker(CH^n(X)\otimes\Z_l\overset{cl}{\To} H^{2n}_\cont(X,\Z_l(n)))\\\to
\Coker(CH^n(X)\{l\}\overset{bl}{\To} H^{2n-1}_\et(X,\Q_l/\Z_l(n)))\to 0
\end{multline*}
o\`u $cl$ est la classe de cycle $l$-adique, $\bar{cl}$ est la classe de
cycle \`a coefficients divisibles et $bl$ est le raffinement de Bloch
de la classe de cycle 
$l$-adique sur les cycles alg\'ebriques de  torsion \cite{bloch1}.
\end{cor}

\prf Cela r\'esulte imm\'ediatement du diagramme commutatif de suites
exactes
\[\begin{CD}
0&&0&&0\\
@V{}V{}V @V{}V{}V @V{}V{}V\\
CH^n(X)\{l\}&=& H^{2n-1}_\Zar(X,\Q_l/\Z_l(n))@>>>
H^{2n-1}_\et(X,\Q_l/\Z_l(n))\\ 
@V{}V{}V @V{}V{}V @V{}V{}V\\
CH^n(X)\otimes\Z_l&=& H^{2n}_\Zar(X,\Z(n))\otimes\Z_l@>>>
H^{2n}_\cont(X,\Z_l(n))\\ 
@V{}V{}V @V{}V{}V @V{}V{}V\\
CH^n(X)\otimes\Q_l&=& H^{2n}_\Zar(X,\Z(n))\otimes\Q_l@>\sim>>
H^{2n}_\cont(X,\Q_l(n))\\ 
@V{}V{}V @V{}V{}V @V{}V{}V\\
CH^n(X)\otimes\Q_l/\Z_l&=&H^{2n}_\Zar(X,\Q_l/\Z_l(n))@>>>
H^{2n}_\et(X,\Q_l/\Z_l(n))\\ 
@V{}V{}V @V{}V{}V\\
0&&0
\end{CD}\]
o\`u, compte tenu  de la
proposition \ref{p4}, les $0$ sup\'erieurs du milieu et de droi\-te
proviennent du corollaire
\ref{c3} c) et e) et l'isomorphisme provient de ce corollaire et du
th\'eor\`eme
\ref{t4}.\qed 

\begin{cor}[\cf\protect{\cite[remark 5.6 a)]{milne2}, \cite[\S
3]{mil-ram}}] Soit $X\in B_\tate(k)$. Pour que l'application cycle
$l$-a\-di\-que enti\`ere
$CH^n(X)\otimes\Z_l\to H^{2n}_\cont(X,\Z_l(n))$ soit surjective, il faut
et il suffit que
\begin{thlist}
\item La classe de cycle \`a coefficients divisibles
\[CH^n(X)\otimes\Q_l/\Z_l\overset{\bar{cl}}{\To}
H^{2n}_\et(X,\Q_l/\Z_l(n))\]
soit injective;
\item L'application de Bloch $CH^n(X)\{l\}\overset{bl}{\To}
H^{2n-1}_\et(X,\Q_l/\Z_l(n))$ soit surjective.\hfill $\square$
\end{thlist}
\end{cor}

\section{Vari\'et\'es ouvertes}\label{ouv}

\begin{lemme}\label{l5} Soit $X$ une $k$-vari\'et\'e projective lisse de
dimension $d$. Alors sous les conditions suivantes
\begin{thlist}
\item $d\le 1$
\item $n\le 1$
\end{thlist}
les groupes
$H^i_W(X,\Z(n))$ sont de type fini pour $i\ne
2n+1$. Ils sont finis pour $i\ne 2n,2n+1$ ainsi que pour $n>d$, et nuls
pour $i\le 0$ si $n>0$.
\end{lemme}

\prf Si $n=0$, l'\'e\-non\-c\'e
r\'e\-sul\-te de \cite[th. 3.2]{licht2}. Supposons
$n=1$: alors
$\Z(1)=\G_m[-1]$. La suite exacte \eqref{e1} se retraduit en une suite
exacte
\begin{multline}\label{e1_1}
\dots\to H^{i-1}_\et(X,\G_m)\to H^i_W(X,\Z(1))\\
\to H^{i-2}_\et(X,\G_m)\otimes\Q\stackrel{\partial}{\to}
H^{i}_\et(X,\G_m)\to\dots
\end{multline}
Le groupe $H^{i-2}_\et(X,\G_m)\otimes\Q$ est nul pour $i\ne 3$. Le
groupe
$H^i_W(X,\Z(1))$ est donc nul pour $i\le 0$, isomorphe \`a $E^*$ pour
$i=1$ (o\`u $E$ est le corps des constantes de $X$) donc fini, et
isomorphe \`a $\Pic(X)$ pour $i=2$, donc de type fini par le th\'eor\`eme
de N\'eron-Severi. Pour $i>3$, il est isomorphe \`a un quotient de
$H^{i-1}_\et(X,\G_m)$; par la suite exacte de Kummer, ce groupe est
lui-m\^eme isomorphe \`a $H^{i-1}_\et(X,(\Q/\Z)'(1))$ qui est fini par
\cite[th. 2]{ctss}. Enfin, si
$X$ est un point, les $H^i_\et(X,\G_m)$ sont finis pour tout $i$, donc
aussi les
$H^i_W(X,\Z(1))$. Ceci d\'emontre (ii).

Si $X$ est un point ou une courbe, toute la strat\'egie de
d\'emonstration conduisant au corollaire \ref{c3} s'applique: la
conjecture de Tate est vraie pour $X$ et l'\'equivalence rationnelle
co\"{\i}ncide avec l'\'equivalence num\'erique, donc le th\'eor\`eme
\ref{t4} est vrai pour $X$. On en d\'eduit (i).\qed

\begin{rem} D'apr\`es \eqref{e1_1}, la g\'en\'eration finie de
$H^3_W(X,\Z(1))$  implique la finitude de $H^2_\et(X,\G_m)=Br(X)$,
elle-m\^eme \'equivalente \`a la validit\'e de la conjecture de Tate en
codimension $1$ pour $X$. Cette implication est donc en fait une
\'equivalence par les raisonnements du paragraphe pr\'ec\'edent.
\end{rem}

\begin{lemme}\label{l6} Soit $X$ une $k$-vari\'et\'e lisse de dimension
$d$. Alors sous les conditions suivantes
\begin{thlist}
\item $d\le 1$
\item $n\le 1$
\end{thlist}
les groupes
$H^i_W(X,\Z[1/p](n))$ sont des $\Z[1/p]$-modules de type fini pour $i\le
n$ ou $i>2n+1$. Ils sont finis pour $i\notin [n,2n+1]$ ainsi que pour
$n>d$, et nuls pour $i\le 0$ si $n>0$.
\end{lemme}

\prf En partant du lemme \ref{l5}, on utilise la m\^eme m\'ethode de
d\'evissage que celle de \cite[d\'em. du th. 1]{lcoh}: c'est loisible
car les groupes $H^i_W(X,\Z[1/p](n))$ satisfont \`a un th\'eor\`eme de
puret\'e, ce qui se r\'eduit \`a la puret\'e de la cohomologie \'etale
\cite[cor. 1.19]{glr} via l'isomorphisme \eqref{e14} (c'est pour avoir
cette puret\'e qu'il est n\'ecessaire d'inverser $p$). Le seul point \`a
v\'erifier est qu'on ne ``perd" pas de g\'en\'eration finie quand on
utilise le th\'eor\`eme de de Jong
\cite{dJ}. 

Soit donc $\tilde U\to U$ un rev\^etement fini de vari\'et\'es lisses,
et supposons que $H^i_W(\tilde U,\Z[1/p](n))$ soit un $\Z[1/p]$-module de
type fini. Par un argument de transfert, $\Ker(H^i_W(U,\Z[1/p](n))\to
H^i_W(\tilde U,\Z[1/p](n)))$ est \hbox{d'exposant} fini, disons $m$. Mais
$H^{i-1}_\et(U,\mu_m^{\otimes n})=H^{i-1}_W(U,\mu_m^{\otimes n})$ est
fini par
\cite[Th. finitude]{deligne2}, ce qui montre que ce noyau est fini, d'o\`u
la g\'en\'eration finie de $H^i_W(U,\Z[1/p](n))$.\qed

La proposition suivante raffine le corollaire \ref{c2}:

\begin{prop}\label{c5} Soit $U$ un ouvert non vide de $X$, o\`u $X\in
B_\tate(k)$, et soit $d=\dim X=\dim U$. Alors, dans les deux cas suivants:
\begin{thlist}
\item $d\le 2$
\item $n\le 2$
\end{thlist}
les conclusions du lemme \ref{l6} sont vraies pour $U$.
\end{prop}

\prf M\^eme m\'ethode que ci-dessus, en utilisant le corollaire \ref{c3}
c) et le lemme \ref{l6}. Plus pr\'ecis\'ement:

Notons $Z=X-U$ (structure r\'eduite) et stratifions $Z$ par
sa cha\^{\i}ne canonique de lieux singuliers ($Z_0=Z,
Z_{r+1}=(Z_r)_{sing}$). Il suffit de voir que, pour tout
$r\ge 0$, le groupe
$H^j_W(Z_r-Z_{r+1},\Z[1/p](n-c_r))$ v\'erifie les hypoth\`eses du
lemme \ref{l6}, o\`u
$c_r$ est la codimension de $Z_r$ dans $X$.

a) Si $d\le 2$, les $Z_r-Z_{r-1}$ sont soit des points soit des courbes,
et on est dans le cas (i) du lemme \ref{l6}.

b) Si $n\le 2$, on a $n-c_r\le 1$ pour tout $r$, et on est dans le
cas (ii) du lemme \ref{l6}.\qed

\begin{thm}\label{c6} Soient $X\in B_\tate(k)$ et
$U$ un ouvert de $X$. Alors les grou\-pes
\[H^0(U,\sK_2),H^1(U,\sK_2),H^2(U,\sK_2)=CH^2(U)\]
sont de type fini. Le groupe $K_3(K)_\ind$ est \'egal \`a $K_3(k)$, o\`u
$K$ est le corps des fonctions de $X$ (ou de $U$).
\end{thm}

\prf Si l'on acceptait d'inverser $p$, ce r\'esultat se
d\'eduirait directement de la proposition \ref{c5} via \cite[th. 1.1 et
1.6]{k}; pour obtenir un r\'esultat exact, on est oblig\'e de refaire la
d\'emonstration (puret\'e) avec la cohomologie de Zariski motivique
plut\^ot qu'avec la cohomologie de Lichtenbaum motivique.

Traitons d'abord le cas de
$X$: d'apr\`es
\cite[th. 1.6]{k}, les groupes cit\'es ne sont autres que $H^i(X,\Z(2))$
pour $i$ respectivement
\'egal \`a $2,3,4$ et $1$. D'autre part, d'apr\`es \loccit, th. 1.1,
l'homomorphisme $H^i(X,\Z(2))\to H^i_\et(X,\Z(2))$ est bijectif pour
$i\le 3$ et injectif pour $i=4$. L'\'enonc\'e pour les $H^i(X,\sK_2)$
r\'esulte donc encore du corollaire \ref{c3} c). De plus, on
obtient que $K_3(K)_\ind$ est de type fini. Mais, d'apr\`es
Merkurjev-Suslin \cite{ms}, l'homomorphisme
\[K_3(k)\to K_3(K)_\ind\]
est injectif \`a conoyau divisible; il est donc bijectif (\cf
\cite{knote}).

(Dans \cite{k}, nous travaillions avec le complexe $\Gamma(2)$ de
Lichtenbaum: les raisonnements sont identiques en rempla\c cant
$\Gamma(2)$ par $\Z(2)$ et en utilisant les r\'esultats de
Merkurjev-Suslin \cite{ms2,ms}.)

Dans le cas g\'en\'eral, on raisonne comme dans la d\'emonstration de la
proposition \ref{c5}, en utilisant le th\'eor\`eme de puret\'e pour la
cohomologie motivique et le fait que les groupes $H^i_\Zar(X,\Z(m))$ sont
de type fini pour $m\le 1$ et toute $k$-vari\'et\'e lisse $X$: pour $m=0$
c'est \'evident et pour $m=1$ cela se r\'eduit \`a la g\'en\'eration
finie de $\Gamma(X,\G_m)$ ($i=1$) et de $\Pic(X)$ ($i=2$).\footnote{On
peut r\'eduire la preuve de la g\'en\'eration finie de $\Pic(X)$ au cas
projectif en passant par le th\'eor\`eme de de Jong, par le m\^eme
raisonnement que dans la d\'emonstration du lemme
\protect\ref{l6}.}\qed

\begin{thm}\label{c4} Soient $X\in B_\tate(k)$ de dimension $\le 2$ et
$K=k(X)$.
Alors on a des isomorphismes canoniques
\[\left(K_n^M(K)\oplus\bigoplus_{0\le i\le n-1}
H^{2i-n-1}(K,(\Q/\Z)'(i))\right)\otimes\Z_{(2)}\iso
K_n(K)\otimes\Z_{(2)}\]  o\`u $(\Q/\Z)'(i)=\colim_{(m,\car k)=1}
\mu_m^{\otimes i}$. De plus,
$K_n^M(K)$ est de torsion pour $n\ge 3$ et nul pour
$n\ge 4$.
\end{thm}

\prf Supposons d'abord $\car k=2$. Alors le second facteur du membre de
gauche est nul. L'isomorphisme $K_n^M(K)\otimes\Z_{(2)}\iso
K_n(K)\otimes\Z_{(2)}$ r\'esulte de \cite{geilevp}.

Supposons maintenant $\car k\ne 2$. Pour construire le morphisme, on part
des isomorphismes de \cite[th. 1]{local}:
\[\bigoplus_{0\le i\le n+1}
H^{2i-n-1}_\et(K,\mu_{2^\nu}^{\otimes i})\overset{\sim}{\To}
K_{n+1}(K,\Z/2^\nu).\]

En prenant la limite inductive sur $\nu$, on obtient des
isomorphismes
\[\bigoplus_{0\le i\le n+1}
H^{2i-n-1}_\et(K,\Q_2/\Z_2(i))\overset{\sim}{\To}
K_{n+1}(K,\Q_2/\Z_2).\]

On obtient l'homomorphisme de l'\'enonc\'e en composant avec le Bockstein 
$K_{n+1}(K,\Q_2/\Z_2)\to K_n(K)\otimes\Z_{(2)}$ et en n\'egligeant les
facteurs $i=n,n+1$.

Pour montrer que cet homomorphisme est un isomorphisme, on raisonne comme
dans \cite{local} en utilisant la suite spectrale (convergente) de
Bloch-Lich\-ten\-baum
\cite{bl}
\begin{equation}\label{e11}
E_2^{p,q}=H^{p-q}(K,\Z(-q))\Rightarrow K_{-p-q}(K).
\end{equation}

\begin{sloppypar}
Tout d'abord, le corollaire \ref{c5} montre que le groupe
$H^i_\Zar(K,\Z(n))$ est de torsion pour $i<n$. Il en r\'esulte que le
Bockstein $H^{i-1}_\Zar(K,(\Q/\Z)'(n))\to H^i_\Zar(K,\Z[1/p](n))$ est un
isomorphisme pour $i<n$. D'autre part, la conjecture de Milnor \cite{voem}
et \cite{geilev} montrent que les homomorphismes
$H^{i-1}_\Zar(K,\Q_2/\Z_2(n))\to H^{i-1}_\et(K,\Q_2/\Z_2(n))$
sont des isomorphismes pour $i\le n$. On en d\'eduit des isomorphismes
\begin{equation}\label{e12}
H^{i-1}_\et(K,\Q_2/\Z_2(n))\overset{\sim}{\To}
H^i_\Zar(K,\Z(n))\otimes\Z_{(2)},\quad i<n.
\end{equation}
\end{sloppypar}

En utilisant la compatibilit\'e de \eqref{e11} aux
produits et aux transferts \cite{fs,levine} et les isomorphismes
\eqref{e12}, et en examinant la construction des homomorphismes du
th\'eor\`eme
\ref{c5}, on d\'emontre comme dans \cite[\S 3]{local} que ces derniers
d\'etruisent successivement les diff\'erentielles de la suite spectrale
\eqref{e11} (localis\'ee en $2$) et scindent la filtration donn\'ee par
celle-ci sur l'aboutissement. Ceci conclut la d\'emonstration, et
donne de plus que la suite spectrale \eqref{e11} d\'eg\'en\`ere
canoniquement apr\`es localisation en $2$.\qed

\begin{cor}\label{c4.2}  L'alg\`ebre
$K_*(K)\otimes\Z_{(2)}$ est engendr\'ee par les unit\'es
et la $K$-th\'eorie du sous-corps des constantes, \`a transfert pr\`es.
\end{cor}

\prf Cela r\'esulte de la construction dans
\cite{local} des fl\`eches donnant naissance \`a l'isomorphisme du
th\'eor\`eme \ref{c4}. Plus pr\'ecis\'ement, l'homomorphisme compos\'e
\begin{multline*}
\bigoplus_{[k':k]<\infty} K_{2i-n-1}^M(k'K)\otimes
H^0(k',\Q_2/\Z_2(n+1-i))\to\\ \bigoplus_{[k':k]<\infty}
H^{2i-n-1}(k'K,\Q_2/\Z_2(i))\overset{\Cores}{\To}
H^{2i-n-1}(k'K,\Q_2/\Z_2(i))
\end{multline*} 
est surjectif pour tout $0\le i\le n-1$ d'apr\`es la conjecture de Milnor
et \cite[th. 1]{k3}. D'autre part, on a d'apr\`es Quillen des
isomorphismes
$K_{2n-2i+1}^M(k')\otimes\Z_{(2)}\iso H^0(k',\Q_2/\Z_2(n+1-i))$.
Enfin, la composition de la surjection induite
\[\bigoplus_{[k':k]<\infty} K_{2i-n-1}^M(k'K)\otimes
K_{2n-2i+1}(k')\otimes\Z_{(2)}\to
H^{2i-n-1}(K,\Q_2/\Z_2(i))\]
avec la fl\`eche du th\'eor\`eme \ref{c4} n'est autre par construction
(\cf \cite{local}) que celle donn\'ee par le produit en $K$-th\'eorie. Le
corollaire en r\'esulte.\qed

\begin{rem}\label{r2} En passant \`a la limite sur les extensions finies
de $k$, on d\'eduit du corollaire \ref{c4.2} que la conjecture de Suslin 
\cite[conj. 4.1 et note]{suslinicm} est vraie apr\`es localisation en
$2$ pour le corps $L=\bar k K$: l'alg\`ebre $K_*(L)\otimes\Z_{(2)}$ est
engendr\'ee par $K_1(L)\otimes\Z_{(2)}$ et $K_*(\bar k)\otimes\Z_{(2)}$.
\end{rem}

\begin{rem} Soit $U$ un ouvert de $X$, o\`u $X\in B_\tate(k)$ et $\dim
X\le 2$. On peut montrer que, pour tout $n\in\Z$
\begin{thlist}
\item $\zeta(U,s)=(1-q^{n-s})^{a_{n}}\varphi(s)$, avec $a_{n}=\sum
(-1)^{r+1}\rg
H_{r}^{c,\et}(U,\Z(n))$  
\item 
$\displaystyle\varphi(n)
=\prod_{r=0}^{2d}\ind(\bar\partial_{r})^{(-1)^ {r+1}}$ \`a une puissance
de $q$ pr\`es
\end{thlist}
o\`u 
\begin{itemize}
\item[a)] $H_{r}^{c,\et}(U,\Z(n))=H^{2d-r}_\et(U,\Z(d-n))$
\item[b)] $\bar\partial_{r}:H_{r+2}^{c,\et}(U,\Z(n))\otimes (\Q/\Z)'\to
H_r^{c,\et}(U,\Z(n))_\tors[1/p]$ est un certain homomorphisme induit par
l'homomorphisme $\partial$ de \eqref{e1}, dont le noyau et le conoyau sont
finis
\item[c)]
$\ind(\bar\partial_{r})=|\Ker\bar\partial_{r}|/|\Coker\bar\partial_{r}|$.
\end{itemize}

Voir \cite[th. 9.16]{tate}, ainsi que le th\'eor\`eme \ref{t7} (iv) et la
d\'emonstration de la proposition \ref{p6}.
\end{rem}

\begin{rem} Si l'on essaye d'\'etendre les r\'esultats du th\'eor\`eme
\ref{c6} aux poids sup\'erieurs, on se heurte d'abord \`a la conjecture
de Bloch-Kato (isomorphisme de la $K$-th\'eorie de Milnor modulo $m$ avec
la cohomologie galoisienne). Pour les besoins de la discussion, supposons
celle-ci connue. Alors le r\'esultat principal de Geisser-Levine
\cite{geilev} implique que le morphisme canonique
\[\Z(n)\to \tau_{\le n+1} R\alpha_*\alpha^*\Z(n)\]
est un quasi-isomorphisme dans $D^{-}((Sm/k)_\Zar)$. Concr\`etement, on en
d\'eduit que, pour toute $k$-vari\'et\'e lisse $X$, l'homomorphisme
\[H^i_\Zar(X,\Z(n))\to H^i_\et(X,\Z(n))\]
est un isomorphisme pour $i\le n+1$ est est injectif pour $i=n+2$. Si $X$
est par exemple un produit de courbes elliptiques, cela implique via le 
corollaire \ref{c3} c) que $H^i_\Zar(X,\Z(n))$ est de type fini pour
$i\le n+2$.

Pour $n\le 2$, ceci couvre toute la cohomologie motivique de $X$: c'est
le th\'eor\`eme \ref{c6} dans le cas projectif. Examinons le cas $n=3$.
On obtient une suite exacte:
\begin{multline*}
0\to H^5_\Zar(X,\Z(3))\to H^5_\et(X,\Z(3))\to
H^0_\Zar(X,R^5\alpha_*\Z(3))\\
\to CH^3(X)\to H^6_\et(X,\Z(3)).
\end{multline*}

Ceci montre que \emph{sous la conjecture de Bloch-Kato et  pour $X\in
B_\tate(k)$, $CH^3(X)$ est de type fini si et seulement si le groupe de
cohomologie non ramifi\'ee
$H^0_\Zar(X,R^5\alpha_*\Z(3))\simeq H^0_\Zar(X,\sH^4_\et((\Q/\Z)'(3))$ est
de type fini}.

Pour rendre la discussion ci-dessus inconditionnelle, on peut localiser
en $l=2$ et utiliser le th\'eor\`eme de Voevodsky \cite{voem}. On peut
aussi localiser en $l=p$ et utiliser le th\'eor\`eme de Geisser-Levine
\cite{geilevp}.

\begin{sloppypar}
Il serait amusant de donner un exemple o\`u la $2$-torsion de
$H^0_\Zar(X,\sH^4_\et((\Q/\Z)'(3))$ est infinie.
\end{sloppypar}
\end{rem}

\section{Retour aux motifs}\label{retour}

\subsection{Motifs purs} Rappelons que $\sA$ d\'esigne la cat\'egorie des
$k$-motifs de Chow
\`a coefficients rationnels, et $\bar \sA$ la cat\'egorie des motifs purs
modulo l'\'equivalence num\'erique (\'egalement \`a coefficients
rationnels). Donnons-nous un nombre premier $l$. Les foncteurs 
\begin{align*}
X&\mapsto H^i_\cont(\bar X,\Q_l(n))\\
X&\mapsto H^i(X,\Q(n))
\end{align*}
se prolongent en des foncteurs
\begin{align*}
M&\mapsto H^i_\cont(\bar M,\Q_l)\\
M&\mapsto H^i(M,\Q)
\end{align*}
par les r\`egles 
\begin{align*}
H^i_\cont(\bar h(X)(j),\Q_l)&=H^{i-2j}_\cont(\bar X,\Q_l(-j))\\
H^i(h(X)(j),\Q)&=H^{i-2j}(X,\Q(-j)).
\end{align*}

\begin{lemme}\label{l5.1} Pour tout $M\in\sA$ et pour tout $i\in\Z$,
l'action de Frobenius sur $H^i_\cont(\bar M,\Q_l)$ est pure de poids $i$;
son polyn\^ome caract\'eristique ne d\'epend pas de $l$.
\end{lemme}

\prf Cela r\'esulte de
l'hypoth\`ese de Riemann \cite{deligne} et de \cite{km}, par r\'eduction
aux vari\'et\'es projectives lisses.\qed

Introduisons maintenant quelques sous-cat\'egories pleines de $\sA$:

\begin{sloppypar}
\begin{itemize}
\item $\sA_{s,l}$: motifs $M$ tels que l'action de
Frobenius sur $H^*(\bar M,\Q_l)$ soit semi-simple.  C'est une
sous-cat\'egorie \'epaisse, donc pseudo-ab\'elienne puisque $\sA$ l'est.
Elle est aussi mono\"{\i}dale et rigide.
\item $\sA_\kim$: motifs de dimension finie au sens de Kimura. Elle est
aussi pseudo-ab\'elienne rigide \cite{ki}.
\item $\sA_{t,l}$: motifs $M$ tels que les homomorphismes
\begin{align*}
\sA(M,\un)\otimes\Q_l&\to H^0(\bar M,\Q_l)^G\\
\sA(M^\vee,\un)\otimes\Q_l&\to H^0(\bar M^\vee,\Q_l)^G
\end{align*}
soient surjectifs, o\`u $G=Gal(\bar k/k)$. Elle est pseudo-ab\'elienne,
stable par dual mais pas a priori par produit tensoriel. 
\end{itemize}
\end{sloppypar}

La cat\'egorie $\sA_{s,l}\cap \sA_\kim$ contient les motifs de
vari\'et\'es ab\'eliennes et leurs tordus \`a la Tate. Conjecturalement,
on a $\sA_{s,l}=\sA_\kim=\sA_{t,l}=\sA$.

\begin{prop}\label{p5} La cat\'egorie $\sA_{st}=\sA_{s,l}\cap\sA_{t,l}$ ne
d\'epend pas de $l$. Pour tout $M\in \sA_{st}$, on a
\[\dim \bar\sA(M,\un) =-\ord_{s=0}\zeta(M,s)\]
o\`u $\zeta(M,s)$ est la fonction z\^eta de $M$ \cite{kleiman}. De plus,
les homomorphismes
\[\sA_{\hom,l}(M,\un)\to \bar\sA(M,\un),\qquad \sA_{\hom,l}(\un,M)\to
\bar\sA(\un,M),\]
o\`u $\sA_{\hom,l}$ d\'esigne la cat\'egorie des motifs purs modulo
l'\'equivalence $\Q_l$-homologique, sont bijectifs.
\end{prop}

\prf Cela r\'esulte de (la d\'emonstration de) \cite[th.
2.9]{Tate}. \qed

\begin{thm}\label{t3} Soit $\sA'=\sA_{st}\cap \sA_\kim$. Cette
sous-cat\'egorie pleine de $\sA$ est pseudo-ab\'elienne et stable par
dual. De plus, pour tout
$M\in \sA'$, les homomorphismes
\[\sA(\un,M)\to \bar\sA(\un,M),\qquad \sA(M,\un)\to \bar\sA(M,\un)\]
sont bijectifs.
\end{thm}

\prf Cela se d\'emontre comme le th\'eor\`eme \ref{t1}.\qed

\begin{prop}\label{c7} Pour $M\in \sA'$, on a $H^i(M,\Q)=H^i_\cont(\bar
M,\Q_l)^G=0$ pour $i\ne 0$ et pour tout $l$.
\end{prop}

\prf Pour la nullit\'e du premier groupe, m\^eme d\'emonstration que
\cite[th. 3.3]{geisser}: par un raisonnement analogue \`a celui de la
d\'emonstration de la proposition
\ref{p1} (utilisant la proposition \ref{p0} et la semi-simplicit\'e de
$\bar \sA$), on peut supposer $M\in\bar \sA$ simple. Si $M\ne \un$, on a
$F_M\ne \un$, d'o\`u l'\'enonc\'e. Si $M=\un$, l'\'enonc\'e est
\'evident. La nullit\'e du deuxi\`eme groupe d\'ecoule du lemme
\ref{l5.1}.\qed

\begin{rem} Soit $\sA_\el$ la sous-cat\'egorie \'epaisse de $\sA$
engendr\'ee par les motifs des produits de courbes elliptiques et leurs
tordus \`a la Tate. Alors $\sA_\el\subset \sA'$ d'apr\`es \cite{spiess}.
Cette sous-cat\'egorie est rigide. Il en r\'esulte via le th\'eor\`eme
\ref{t3} (ou par le th\'eor\`eme \ref{t1}) que le foncteur induit
\[\sA_\el\to \bar\sA_\el\]
est une \'equivalence de cat\'egories. En particulier, $\sA_\el$ est
ab\'elienne semi-simple.
\end{rem}

Nous allons maintenant donner un analogue du corollaire  \ref{c3l} pour
les objets de $\sA'$. Pour cela, nous avons besoin d'introduire une
autre cat\'egorie: les \emph{motifs de Chow \'etales}.

\begin{defn}\label{chowet} La cat\'egorie $\Chow_\et$ des motifs de
Chow \'etales est d\'efinie de la mani\`ere suivante:
\begin{itemize}
\item[a)] On introduit la cat\'egorie des correspondances de Chow \'etales
effectives
$\Cores_\et^\eff$, dont les objets sont
$k$-vari\'et\'es projectives lisses ($X\mapsto [X]$) et les morphismes des
\'el\'ements des groupes
\[\Cores_\et([X],[Y])=H^{2\dim X}_\et(X\times Y,\Z(\dim X))\]
la composition \'etant induite par la formule habituelle.
\item[b)] On construit la cat\'egorie des motifs de Chow \'etales
effectifs
$\Chow_\et^\eff$: un objet est un couple $([X],p)$, o\`u $p=p^2\in
\Cores_\et^\eff([X],[X])\otimes\Q$; un morphisme de $([X],p)$ vers
$([Y],q)$ est un \'el\'ement $f\in \Cores_\et^\eff([X],[Y])$ tel que
$\hat{f}=q\hat{f}p$, o\`u $\hat{f}$ est l'image de $f$ dans
$\Cores_\et^\eff([X],[Y])\otimes\Q$; la composition des morphismes est
induite par celle de $\Cores_\et^\eff$.
\item[c)] On passe de $\Chow^\eff_\et$ \`a $\Chow_\et$ en inversant le
motif de Lefschetz. On note $h:SmProj/k\to \Chow_\et$ le foncteur
\'evident.
\end{itemize}
\end{defn}

La justification que toutes ces constructions est
donn\'ee dans l'appendice \ref{A}. La cat\'egorie $\Chow_\et$ est
pseudo-ab\'elienne, mono\"{\i}dale sy\-m\'e\-tri\-que ri\-gi\-de. On a un
foncteur \'evident (additif, mono\"{\i}dal)
\[\sA\to (\Chow_\et)\otimes\Q\]
qui est \emph{pleinement fid\`ele} par la proposition \ref{p4} et m\^eme
\emph{essentiellement surjectif} puisque les deux cat\'egories sont
pseudo-ab\'eliennes. On en d\'eduit un foncteur (additif, tensoriel)
\begin{equation}\label{e7}
\Chow_\et\to\sA
\end{equation}
bien d\'efini \`a isomorphisme naturel mono\"{\i}dal pr\`es, et
\emph{essentiellement surjectif} par construction.

On notera que si $([X],p)\in \Cores^\eff_\et$, il existe $m\in\Z-\{0\}$ tel
que $mp$ provienne d'un $\tilde p\in  \Cores^\eff_\et([X],[X])$. Un tel
$\tilde p$ induit des morphismes $[X]\to ([X],p)$ et $([X],p)\to [X]$,
dont la composition dans les deux sens est (en choisissant $m$ assez
grand) la multiplication par $m$. Par contre, $([X],p)$ n'est pas en
g\'en\'eral facteur direct de $[X]$.

\begin{prop}\label{p5.8} Le foncteur bigradu\'e $X\mapsto H^i_W(X,\Z(n))$
des vari\'et\'es projectives lisses vers les groupes ab\'eliens induit un
foncteur gradu\'e
\begin{align*}
\Chow_\et&\to Ab^*\\
M&\mapsto H^*_W(M,\Z)
\end{align*}
tel que $H^i_W(h(X)(n),\Z)=H^{i-2n}_W(X,\Z(-n))$
pour toute vari\'et\'e projective lisse $X$ et tout $n\in \Z$. \end{prop}

\prf Voir \ref{A.11}. Plus pr\'ecis\'ement, si
$([X],p)\in \Cores^\eff_\et$, on d\'efinit
\[H^i_W(([X],p),\Z(n))= \{x\in H^i_W(X,\Z(n))\mid p\hat{x}=\hat{x}\}\]
o\`u $\hat{x}$ est l'image de $x$ dans $H^i_W(X,\Z(n))\otimes\Q$.\qed

De m\^eme, pour tout $l$, la cohomologie \'etale continue induit un
foncteur gradu\'e
\begin{align*}
\Chow_\et&\to {\Z_l-Mod}^*\\
M&\mapsto H^*_\cont(M,\Z_l)
\end{align*}
et on a une transformation naturelle
\begin{equation}\label{e15}
H^*_W(-,\Z)\otimes\Z_l\to H^*_\cont(-,\Z_l).
\end{equation}

Afin d'\'enoncer le th\'eor\`eme suivant, d\'efinissons des groupes
$H^j(M,\Omega^i)$ pour tout $M\in \Chow_\et$: d'une part, les
correspondances de Chow \'etales op\`erent sur la cohomologie de Hodge
puisque celle-ci peut s'exprimer en termes de cohomologie \'etale.
D'autre part, pour toute vari\'et\'e projective lisse $X$, on a un
isomorphisme
\[H^j(X,\Omega^i)\oplus
H^{j-1}(X,\Omega^{i-1})\overset{\sim}{\To} H^j(X\times\P^1,\Omega^i)\]
induit par le cup-produit par la ``premi\`ere classe de Chern" de
$\sO(1)$ dans $H^1(\P^1,\Omega^1)$. Ceci justifie le fait que la
d\'efinition
\[H^j(h(X)(n),\Omega^i)=H^{j-n}(X,\Omega^{i-n})\]
s'\'etend fonctoriellement \`a $\Chow_\et$ par passage aux facteurs
directs.

\begin{thm}\label{t7}
Soit $\Chow'_\et$ l'image r\'eciproque pleine de $\sA'$ par le foncteur
\eqref{e7}. Alors,\\
a) Pour tout $l$, la restriction de \eqref{e15} \`a $\Chow'_\et$ est
un isomorphisme de foncteurs.\\
b) Pour tout $M\in\Chow'_\et$, les groupes $H^i_W(M,\Z)$
sont de type fini. Ils sont finis pour $i\ne 0,1$ et
nuls sauf pour un nombre fini de valeurs de $i$. L'homomorphisme 
$H^0_W(M,\Z)\overset{\cdot e}{\To} H^1_W(M,\Z)$ a un noyau et un conoyau
finis. De plus:
\begin{thlist}
\item L'accouplement $H^0_W(M,\Z)\times H^0_W(M^\vee,\Z)\to \Z$
induit par la fl\`eche de dualit\'e $\un\to M\otimes M^\vee$ est parfait
modulo torsion.
\item Les accouplements analogues
\[H^{i+1}_W(M,\Z)_\tors\times H^{-i}_W(M^\vee,\Q/\Z)\to\Q/\Z\] 
induisent des accouplements parfaits de groupes finis
\[H^{i+1}_W(M,\Z)_\tors\times H^{-i+1}_W(M^\vee,\Z)_\tors\to\Q/\Z.\]
\item L'homomorphisme canonique $\Chow'_\et(M,\un)\to H^0_W(M,\Z)$
est bijectif.
\item Soit $\hat M$ l'image de $M$ dans
$\sA'$. Alors l'ordre du p\^ole de
$\zeta(\hat M,s)$ est \'egal \`a $m=\rg H^0_W(M,\Z)=\dim \sA'(\hat
M,\un)$. On a
\[\lim_{s\to 0} (1-q^{-s})^m\zeta(\hat M,s)=\pm q^{\chi(M,\sO)}\prod_i
|H^i_W(M,\Z)_\tors|^{(-1)^i}\cdot R(M)^{-1}\]
o\`u $R(M)$ est la valeur absolue du d\'eterminant de l'accouplement de
(i) modulo torsion par rapport \`a des bases quelconques des groupes
accoupl\'es et
\[\chi(M,\sO)=\sum_{i,j} (-1)^{i+j+1}i{h_{ij}}, h_{i,j}=\dim_k
H^j(M,\Omega^i).\]
En particulier, le second membre ne d\'epend que de $\hat M$.
\end{thlist}
\end{thm}

\prf On proc\`ede comme au \S 3.\qed

\begin{rems}[concernant le th\'eor\`eme \protect{\ref{t7}} (iv)]\
\begin{enumerate}
\item Ceci donne un sens
pr\'ecis \`a \cite[conj. 7.3]{milne3}, o\`u le terme $\chi(M,\sO)$
(not\'e $\chi(M,\sO,0)$) n'\'etait pas d\'efini. En fait, les groupes
$H^j(\hat{M},
\Omega^i)$ n'ont pas de sens individuellement pour $\hat{M}\in \sA$,
puisqu'ils sont annul\'es par $p$. Il n'est pas clair que $\chi(M,\sO)$
ne d\'epende que de $\hat{M}$.
\item Pour $n\in\Z$, la valeur sp\'eciale de $\zeta(\hat M,s)$ en $s=n$
s'obtient en appliquant ce r\'esultat \`a $M(-n)$.
\item On peut pr\'eciser le signe dans l'expression donn\'ee, exactement
com\-me dans la remarque \ref{r1}.
\item Le foncteur $\Chow_\et\to\sA$ \'etant essentiellement surjectif, le
th\'eor\`eme \protect{\ref{t7}} (iv) d\'ecrit les valeurs
sp\'eciales de $\zeta(\hat M,s)$ pour tous les motifs $\hat M\in\sA'$.
\end{enumerate}
\end{rems}

Pour $M\in \Chow'_\et$, posons $H^q_W(M,\Z(n))=H^{q-2n}_W(M(-n),\Z)$.

\begin{sloppypar}
\begin{cor}[``principe d'identit\'e de Manin arithm\'etique"]
\label{idmanar} Soit
$f:M\allowbreak\to M'$ un morphisme de
$\Chow'_\et$. Supposons que
$f^*:H^q_W(M',\Z(n))\to H^q_W(M,\Z(n))$ et
$f_*:H^q_W(M^\vee,\Z(n))\to H^q_W((M')^\vee,\Z(n))$ soient des
isomorphismes pour tout $q,n\in\Z$ (de mani\`ere \'equivalente via
\eqref{e1}, remplacer la cohomologie de Lichtenbaum par la cohomologie
\'etale). Alors\\ 
a) $\zeta(\hat M,s)=\zeta(\hat M',s)$;\\ 
b) $f$ induit un isomorphisme dans $\sA'$. (On pourrait dire que $f$ est
une \emph{isog\'enie}.)
\end{cor}
\end{sloppypar}

(Ce r\'esultat est \emph{faux} si l'on prend  des coefficients rationnels:
prendre $M=0$ et pour $M'$ un motif simple qui n'est pas une puissance du
motif de Lefschetz.)

\prf Soit $f\in \Q(t)$ tel que $f(n)$ soit d\'efini et \'egal \`a
 $\pm 1$ pour tout $n\ge 0$: alors $f$
est constant, \'egal \`a $\pm 1$. Ceci montre avec le th\'eor\`eme
\ref{t7} b) (iv) que, sous les hypoth\`eses,
$\zeta(M,s)/\zeta(M',s)=\pm 1$. Mais
$\lim\limits_{s\to+\infty}\zeta(M,s)=\lim\limits_{s\to+\infty} 
\zeta(M',s)=1$, d'o\`u a). En utilisant l'hypoth\`ese de Riemann, il
s'ensuit que $f^*$ induit des isomorphismes sur tous les  groupes de
cohomologie $\Q_l$-adique g\'eom\'etriques, et donc sur les $H^i(-,\Q)$
d'apr\`es le th\'eor\`eme \ref{t7} a) et b) (iii), ainsi que la
proposition \ref{p4}. b) r\'esulte maintenant du principe d'identit\'e de
Manin (g\'eom\'etrique).\hfill $\square$

\begin{qn} Est-ce que $f$ est m\^eme un isomorphisme dans $\Chow'_\et$?
\end{qn}

\subsection{Motifs mixtes} Notons $\sD$ la cat\'egorie
$DM_\gm(k)\otimes\Q$ de Voevodsky \cite{voetri}. Par
ailleurs, choisissons une sous-cat\'egorie \'epaisse
$\sA''\subset\sA'$ qui soit  \emph{rigide}, c'est-\`a-dire stable par
produit tensoriel: on peut par exemple prendre $\sA''=\sA_\el$, ou pour
$\sA''$ la sous-cat\'egorie additive tensorielle \'epaisse engendr\'ee par
le motif d'une vari\'et\'e ab\'elienne simple ``de type K3" ou  ``presque
ordinaire" (voir les exemples de l'introduction).
D'apr\`es le th\'eor\`eme
\ref{t3},
$\sA''$ est ab\'elienne semi-simple. Soit
$\sD''$ la sous-cat\'egorie
\'epaisse (triangul\'ee) de $\sD$ engendr\'ee par l'image essentielle de
$\sA''$ via le foncteur de \loccit, prop. 2.1.4. Notons $\delta: \sA\to
\sD$ ce foncteur, ainsi que sa restriction \`a $\sA''$.

\begin{prop} \label{propmot} Pour tous $M,N\in \sA''$, on a
\[\sD''(\delta(M),\delta(N)[i])=
\begin{cases}
0&\text{pour $i\ne 0$}\\
\sA''(M,N)&\text{pour $i=0$.}
\end{cases}\]
\end{prop}

\prf Vu le corollaire \ref{c7}, il suffit de d\'emontrer
que pour $M,N\in\sA$, on a un isomorphisme canonique
\[\sD'(\delta(M),\delta(N)[i])=H^i(M\otimes N^\vee,\Q).\]

Par dualit\'e, on se ram\`ene au cas $N=\un$. On se ram\`ene ensuite au
cas o\`u $M$ est de la forme $h(X)(n)$ o\`u $X$ est projective lisse. La
formule r\'esulte alors de \cite[cor. 3.2.7]{voetri} et du
th\'eor\`eme de simplification \cite{voecan}, qui est maintenant valable
sur un corps parfait quelconque.\qed

Le th\'eor\`eme suivant est une variation sur le th\`eme de
\cite[th. 2.49]{milne}.

\begin{thm}\label{zetamot} La cat\'egorie triangul\'ee $\sD''$ est
semi-simple: tout triangle exact est somme directe de triangles scind\'es.
Le foncteur $\delta$
induit une \'equivalence de cat\'egories
\begin{align*}
\Delta:(\sA'')^{(\N)}&\to \sD''\\
(M_i)&\mapsto \bigoplus_i \delta(M_i)[i].
\end{align*}
\end{thm}

\prf En quatre \'etapes:

1) $\Delta$ est pleinement fid\`ele. En effet, pour 
$(M_i),(N_i)\in (\sA'')^{(\N)}$, on a
\[(\sA'')^{(\N)}((M_i),(N_i))=\bigoplus_{i} (\sA'')(M_i,N_i)\]
et
\begin{multline*}\sD''(\bigoplus_i \delta(M_i)[i],\bigoplus_j 
\delta(N_j)[j])=\bigoplus_{i,j}\sD''(\delta(M_i),\delta(N_j)[j-i])\\ 
=\bigoplus_i\sD''(\delta(M_i),\delta(N_i))
=\bigoplus_i\sA''(M_i,N_i)\end{multline*}
par la proposition \ref{propmot}.

2) L'image essentielle de $\Delta$ est stable par triangles. Cela
r\'esulte facilement du fait que dans une cat\'egorie ab\'elienne
semi-simple, tout morphisme est somme directe d'un morphisme nul et d'un
isomorphisme (\cf par exemple \cite[lemme A.2.13]{ak}).

3) $\Delta$ est essentiellement surjectif. Cela d\'ecoule de 2) et
de la d\'efinition de $\sD''$.

4) Il reste \`a voir que $\sD''$ est semi-simple: cela d\'ecoule encore de
\cite[lemme A.2.13]{ak}.\qed

\begin{cor} Le foncteur $\delta$ 
induit une \'equivalence de cat\'egories
\[D^b(\sA'')\iso \sD''.\]
\end{cor}

\prf En effet, puisque $\sA''$ est semi-simple, 
$D^b(\sA'')$ co\"{\i}ncide avec le membre de gauche de 
l'\'equivalence de cat\'egories $\Delta$ du th\'eor\`eme \ref{zetamot}.
\qed

\begin{cor}\label{c8} Soit $U\in Sm/k$ tel que $M(U)\in\sD''$. Alors, pour
tout
$l\ne p$, l'action de Frobenius sur $H^*_{c}(\overline U,\Q_l)$ est
semi-simple. De plus, le polyn\^ome caract\'eristique de
cette action est ind\'ependant de $l$.
\end{cor}

\prf Soit $d=\dim U$. Par le th\'eor\`eme \ref{zetamot}, on peut \'ecrire
$M^c(U)\allowbreak :=M(U)^*(d)[2d]$ comme somme directe de motifs de la
forme
$\delta(M)[i]$, avec
$M\in\sA''$. Le corollaire r\'esulte alors du lemme \ref{l5.1} (pour
l'ind\'ependance de $l$) et du fait que $\sA''\subset \sA_{s,l}$ (pour la
semi-simplicit\'e).\qed

Nous ne sommes pas en mesure d'appliquer directement le corollaire
\ref{c8} \`a un ouvert d'une surface ab\'elienne. Toutefois, on peut lui
appliquer la m\^eme m\'ethode, au prix d'un petit effort.

\begin{prop}\label{p6} Soit $X\in B_\tate(k)$, avec $d=\dim X\le 2$. Alors
tout ouvert dense $U$ de $X$ v\'erifie la conclusion
du corollaire \ref{c8}.
\end{prop}

\prf Le cas $d\le 1$ est facile et laiss\'e au lecteur. Supposant $d=2$,
montrons que $M^c(U)$ est somme directe de motifs de la forme
$\delta(M)[i]$, avec $M\in\sA'$. Par r\'ecurrence n\oe th\'erienne, il
suffit de montrer que, si $U=U'-Z$ avec
$U'$ ouvert et $Z$ ferm\'e lisse, l'\'enonc\'e pour $U'$ implique
l'\'enonc\'e pour
$U$. Notons donc $M(U')=\bigoplus \delta(M_\alpha)[i_\alpha]$,
$M_\alpha\in \sA'$.

Supposons d'abord $\dim Z=1$. Consid\'erons
la compl\'et\'ee projective lisse $C$ de $Z$: on a
\[\sD(M(C),M(U'))=\bigoplus \sD(M(C),\delta(M_\alpha)[i_\alpha]).\]

En utilisant le fait que la conjecture de Tate est vraie
pour les produits de $3$ courbes \cite[th. 3]{soule}, le
raisonnement de la d\'emonstration du corollaire \ref{c2} montre que
$\sD(M(C),\delta(M_\alpha)[i_\alpha])=0$ pour $i_\alpha\ne 0$. Pour
$i_\alpha=0$, on obtient en utilisant \cite{voecan} des isomorphismes
\[\bar\sA(\bar
h(C),M_\alpha)\osi\sA(h(C),M_\alpha)\iso\sD(M(C),\delta(M_\alpha)).\]

Par semi-simplicit\'e de $\bar\sA$, tout morphisme de $M(C)$ vers
$\delta(M_\alpha)$ est donc somme directe d'un morphisme nul et d'un
isomorphisme. On en d\'eduit le pas de r\'ecurrence par d\'evissage
(utilisant le triangle exact de Gysin de \cite[prop. 3.5.4]{voetri}). Le
cas o\`u $\dim Z=0$ est semblable et plus facile.\hfill $\square$

\section{Prospective}\label{pros}

Rappelons le r\'esultat principal de \cite{glr}:

\begin{thm}[\protect{\cite[th. 3.4]{glr}}]\label{t5} Les trois conjectures
suivantes sont \'e\-qui\-va\-len\-tes:
\begin{thlist}
\item Le th\'eor\`eme \ref{t1} (avec $X'=\Spec k$) et le corollaire
\ref{c1} sont vrais  pour toute $k$-vari\'et\'e projective lisse $X$.
\item Pour tout $n> 0$, le morphisme $\Phi_n$ du lemme \ref{l2}
(relatif \`a un nombre premier $l\ne p$ donn\'e) est un
quasi-isomorphisme.
\item Pour tout $n> 0$, le complexe $\Z_l(n)^c_\et$ (relatif \`a un
nombre premier $l\ne p$ donn\'e) est mall\'eable (voir \cite[d\'ef.
2.16]{glr} pour la d\'efinition de mall\'eable).
\end{thlist}
\end{thm}

Les m\'ethodes pr\'ec\'edentes conduisent en fait \`a d'autres
formulations de ces conjectures:

\begin{sloppypar}
\begin{thm}\label{t6} Les conjectures du th\'eor\`eme \ref{t5} sont
encore \'equivalentes aux suivantes:
\begin{thlist}
\item[\rm (iv)] Pour tout $n> 0$, le morphisme \eqref{e3} (relatif \`a un
nombre premier $l\ne p$ donn\'e) est un quasi-isomorphisme dans $\bar
D(\sT_{\Z})$ (\cf mise en garde \ref{meg1}).
\item[\rm (iv bis)] Pour tout $n> 0$, le morphisme \eqref{e4} ou
\eqref{e4p} (relatif
\`a un nombre premier $l$ donn\'e) est un isomorphisme pour toute
$k$-vari\'et\'e projective lisse $X$.
\item[\rm (iv ter)] Pour tout $n> 0$, tout nombre premier $l$ et toute
$k$-vari\'et\'e projective lisse $X$,
les morphismes \eqref{e4} et \eqref{e4p} sont des isomorphismes.
\item[\rm (v)] Pour toute $k$-vari\'et\'e projective lisse $X$, les
groupes $H^i_W(X,\Z(n))$ sont de type fini.
\item[\rm (v bis)] Pour toute $k$-vari\'et\'e lisse $X$, les
groupes $H^i_W(X,\Z[1/p](n))$ sont des $\Z[1/p]$-modules de type fini.
\item[\rm (v ter)] Pour toute $k$-vari\'et\'e projective lisse $X$, les
groupes
$H^i_W(X,\Z_{(l)}(n))$ sont des $\Z_{(l)}$-modules de type fini, o\`u $l$
est un nombre premier donn\'e.
\item[\rm (vi)] Pour toute $k$-vari\'et\'e projective lisse $X$, la
forme forte  de la	conjecture de Tate (concernant $\zeta(X,s)$) est vraie
et le motif de Chow
$h(X)$ est de dimension finie au sens de Kimura\footnote{Kimura
conjecture le second
\'enonc\'e sur tout corps de base: cela r\'esulterait des conjectures
standard et de l'existence d'une d\'ecomposition de Chow-K\"unneth \`a la
Murre, \cf
\protect\cite[ex. 9.2.4]{ak}}. (Avec les notations du th\'eor\`eme
\ref{t3}: $\sA'=\sA$.)
\item[\rm (vi bis)] Pour toute $k$-vari\'et\'e projective lisse $X$, la
forme forte  de la	conjecture de Tate est vraie
et l'alg\`ebre des correspondances de Chow $A=CH^{\dim X}(X\times
X)\otimes\Q$ est \emph{semi-primaire}: son radical de Jacobson $R$ est
nilpotent et $A/R$ est semi-simple.
\item[\rm (vi ter)] La conjecture de
Tate (cohomologique, relative \`a un nombre premier $l$ donn\'e) est
vraie pour les $k$-vari\'et\'es ab\'eliennes et la cat\'egorie rigide
$\sA$ des
$k$-motifs de Chow \`a coefficients rationnels est engendr\'ee par les
motifs de vari\'et\'es ab\'eliennes et les motifs d'Artin.
\item[\rm (vi quater)] La forme forte de la conjecture de Tate est vraie
pour toutes les vari\'et\'es projectives lisses et il n'existe pas de
``motif de Chow fant\^ome": Si $M\in \sA$ est tel que $\bar M=0$ dans
$\bar
\sA$, alors $M=0$.
\end{thlist}
\end{thm}
\end{sloppypar}

\prf (ii) $\iff$ (iv): cela r\'esulte du lemme \ref{l2}, puisque le
foncteur $\bar R\gamma_*$ est conservatif.

(i) $\If$ (iv ter): on proc\`ede comme dans la preuve du th\'eor\`eme
\ref{t4}.

(iv ter) $\If$ (v): on proc\`ede comme dans la preuve du corollaire
\ref{c3}.

(v) $\If$ (v bis): on proc\`ede comme dans la d\'emonstration du lemme
\ref{l6}.

(v bis) $\If$ (v ter) pour $l\ne p$: c'est \'evident.

(v) $\If$ (v ter): c'est \'evident.

(v ter) $\If$ (iv bis): notons $K(n)$ le c\^one de \eqref{e3}.
L'hypoth\`ese implique que, pour toute $k$-vari\'et\'e projective lisse
$X$, les groupes $H^i_W(X,K(n))$ sont des $\Z_l$-modules de type fini.
D'autre part, on sait que \eqref{e3}$\oo\limits^L\Z/l$ est un
quasi-isomorphisme (\cite[th. 1.5]{geilev} ou \cite{geilevp} selon que
$l\ne p$ ou que $l=p$). Par cons\'equent, les
$H^i_W(X,K(n))$ sont uniquement divisibles, donc nuls.

(iv bis) $\If$ (iv): on proc\`ede comme dans \cite[lemme 3.8]{glr}:
r\'eduction \`a \eqref{e3}$\otimes\Q$, puis puret\'e et th\'eor\`eme de
de Jong.

(i) $\If$ (vi ter): la premi\`ere partie de (vi ter) r\'esulte de (i)
d'apr\`es
\cite{Tate}. La deuxi\`eme partie est une cons\'equence classique de la
conjecture de Tate forte pour toutes les vari\'et\'es (qui
r\'esulte de (i)) pour $\bar\sA$, \cf \cite[rem. 2.7]{milne}; mais (i)
implique que $\sA\to\bar\sA$ est une \'equivalence de cat\'egories.

(vi ter) $\If$ (vi): c'est clair, puisque le motif d'une vari\'et\'e
ab\'elienne est de dimension finie au sens de Kimura \cite[ex. 9.1]{ki},
ainsi que tout motif d'Artin, que la forme forte de la conjecture de Tate
r\'esulte de la forme cohomologique et de la semi-simplicit\'e de
l'action de Frobenius sur la cohomologie \cite{Tate}, et que cette
derni\`ere est vraie pour les vari\'et\'es ab\'eliennes \cite{weil}.

(vi) $\If$ (vi bis): cela r\'esulte de \cite[prop. 9.1.14]{ak} et du
th\'eor\`eme de Jannsen \cite{jannsen}.

(vi bis) $\If$ (i): on proc\`ede comme dans la d\'emonstration du
th\'eor\`eme \ref{t1} et du corollaire \ref{c1}.

(i) $\If$ (vi quater) est clair. (vi quater) $\If$ (vi ter): notons
$\sA_\ab$ la sous-cat\'egorie rigide de $\sA$ engendr\'ee par les motifs
de vari\'et\'es ab\'eliennes et les motifs d'Artin, et $\bar \sA_\ab$
son image dans $\bar\sA$. La forme forte de la conjecture de Tate implique
que $\bar\sA_\ab=\bar \sA$ (voir ci-dessus). Soit $X$ une vari\'et\'e
projective lisse: il existe des isomorphismes inverses $\bar f:\bar
h(X)\iso \bar N$, $\bar g:N\iso \bar h(X)$, avec $N\in
\sA_\ab$. Relevons $\bar f$ et $\bar g$ en des morphismes $f,g$ de $\sA$.
Comme $N$ est de dimension finie au sens de Kimura, l'alg\`ebre $\sA(N,N)$
est semi-primaire, donc $1_N-fg$ est nilpotent et quitte \`a modifier
$f$ ou $g$ on peut supposer que $fg=1_N$. On a alors $h(X)\simeq
N\oplus M$, avec $M$ fant\^ome.  
\qed

La formulation de la conjecture des th\'eor\`emes \ref{t5} et
\ref{t6} qui nous semble la plus prometteuse est (vi ter), mais nous nous
garderons bien de faire une conjecture sur la d\'e\-mons\-tra\-tion d'une
conjecture!

\appendix

\section{Fonctorialit\'e de la cohomologie motivique
\'etale}\label{A}

Soient $k$ un corps commutatif et $\sV$ la cat\'egorie des
$k$-vari\'et\'es projectives lisses. Pour qu'une ``th\'eorie
cohomologique" $X\mapsto H(X)$ permette de d\'efinir une cat\'egorie
additive mono\"{\i}dale sym\'etrique de correspondances $\sC$ \`a
partir de $\sV$, avec $\sC(X,Y)=H(X\times_k Y)$, il suffit (et il faut,
essentiellement) qu'elle poss\`ede la fonctorialit\'e suivante:
\begin{itemize}
\item[a)] Images inverses pour les morphismes de projection
\item[b)] Images directes pour les morphismes de projection
\item[c)] Cup-produits
\end{itemize}
soumise aux relations suivantes:
\begin{itemize}
\item[(i)] a) est fonctoriel
\item[(ii)] b) est fonctoriel
\item[(iii)] c) est associatif et commutatif
\item[(iv)] a) commute \`a c)
\item[(v)] Condition d'\'echange pour a) et b) 
\item[(vi)] a), b) et c) interagissent via la formule de projection.
\end{itemize}

Pour que $H$ permette de plus de d\'efinir une cat\'egorie de motifs
(rigide) $\sM$, il suffit qu'elle v\'erifie une ``formule du fibr\'e
projectif" pour $H^*(\P^1)$.

Enfin, si $U:X\mapsto U(X)$ est une ``th\'eorie cohomologique" sur $\sV$
telle que $H\oplus U$ poss\`ede la fonctorialit\'e ci-dessus, le
cup-produit \'etant identiquement nul sur la composante $U\otimes U$
(c'est-\`a-dire que $U$ est un ``module" sur $H$), alors $U$ s'\'etend en
un foncteur additif sur $\sM$.

\begin{sloppypar}
Lorsque $H(X)=CH(X)$ (anneau de Chow), la fonctorialit\'e
ci-dessus r\'esulte de la th\'eorie de l'intersection. Pour donner un
sens \`a la cat\'egorie $\Chow_\et(k)$ de la d\'efinition \ref{chowet},
il faut la justifier pour $H(X)=\bigoplus H^{2n}_\et(X,\Z(n))$: c'est
l'objet de cet appendice. Pour cela, nous justifierons cette
fonctorialit\'e pour tous les ``groupes de Chow
sup\'erieurs \'etales $H^{i}_\et(X,\Z(n))$": en effet nous n'avons pas
trouv\'e de m\'ethode plus simple que de traiter ce cas g\'en\'eral.
\end{sloppypar}

On note $\alpha$ la projection du gros site \'etale de $\Spec k$ sur son
gros site Zariski, et $\alpha_U$ la restriction de $\alpha$ au petit site
\'etale d'un $k$-sch\'ema $U$. Comme d'habitude, on pose
$\Z(n)_U^\et=\alpha^*\Z(n)_U$ pour $n\ge 0$ et
$\Z(n)^\et_U=(\Q/\Z)'(n)_U[-1]$ pour $n<0$, de sorte qu'on a des
morphismes canoniques
\begin{equation}\label{eqA.1}
(\Q/\Z)'(n)_U\to \Z(n)^\et_U[1]
\end{equation}
pour tout $U\in (\Spec k)_\Et$ et tout $n\in\Z$.

\subsection{Fonctorialit\'e des groupes de Chow sup\'erieurs} Nous nous
reposerons sur la fonctorialit\'e correspondante pour les groupes
$H^i_\Zar(X,\Z(n))$. Cette fonctorialit\'e n'est que partiellement
r\'edig\'ee:
\begin{itemize}
\item a) est construit pour tout morphisme vers une $k$-vari\'et\'e lisse
dans \cite[partie I, ch II, 3.5]{mlevine}, et ind\'ependamment par Bloch
(non publi\'e). Ceci corrige une d\'emonstration incompl\`ete de
\cite{bloch}.
\item b) est \'el\'ementaire et existe pour tout morphisme propre
\cite{bloch}.
\item c) est construit sur toute vari\'et\'e lisse dans \cite[\S
8]{geilev}, corrigeant la construction incompl\`ete de \cite{bloch}.
\item (i) est v\'erifi\'e pour des morphismes entre vari\'et\'es lisses
dans \cite[1.2]{glr} (\cf th. 1.17 de \loccit).
\item (ii) est \'el\'ementaire modulo les relations habituelles de
multiplicit\'es.
\item (iii) est bri\`evement mentionn\'e \`a la fin de \cite[\S
8]{geilev}.
\item (iv), (v) et (vi) ne sont pas r\'edig\'es. Marc Levine m'a
inform\'e qu'ils seront incorpor\'es \`a \cite{marc}, ainsi qu'une
version d\'etaill\'ee de (iii).
\end{itemize}

Dans ce qui suit, \emph{nous admettrons les compatibilit\'es (iii), (iv),
(v) et (vi) ci-dessus}. La situation est la suivante: \`a partir du cas
Zariski, la construction de a) et c) en cohomologie motivique \'etale est
``formelle", ainsi que la v\'erification de (i) et (iii). En revanche, la
construction de b) n'est pas formelle: nous la r\'ealiserons pour tout
morphisme projectif par r\'eduction aux immersions ferm\'ees et aux
projections \`a partir d'un espace projectif. 

\subsection{ a): images inverses} La d\'efinition en cohomologie de
Zariski est obtenue \`a partir de zig-zags finis entre complexes
explicites. Il suffit de v\'erifier que ces constructions commutent \`a
la localisation \'etale, ce qui est imm\'ediat.

\subsection{ i): fonctorialit\'e des images inverses} M\^eme
justification que ci-dessus.

\subsection{ c): produits} M\^eme justification (contempler par exemple le
zig-zag de \cite[(8.1)]{geilev} pour deux $k$-sch\'emas $X,Y$
essentiellement de type fini et deux entiers $q,q'$:
\begin{multline*} z^q(X,*)\otimes z^{q'}(Y,*)\to \Tot z^{q+q'}(X\times
Y,*,*)\leftarrow
\Tot z^{q+q'}(X\times Y,*,*)_\sT\\
\to z^{q+q'}(X\times Y,*).)
\end{multline*}

\subsection{ iii): associativit\'e, commutativit\'e et iv):
compatibilit\'e de a) et c)}\label{A.5} M\^eme justification.

\subsection{Formule du fibr\'e projectif}\label{fibproj} Soient $E$ un
fibr\'e vectoriel sur une $k$-vari\'et\'e lisse $Y$ et
$p:\P(E)\to Y$ la projection. En proc\'edant comme dans
\cite[\S 5]{k}, on construit des isomorphismes
\[\bigoplus_{j=0}^{N-1} \Z(n-j)^\et_Y[-2j]\iso
Rp_*\Z(n)^\et_{\P(E)}\]
avec $N=\dim E$. Plus pr\'ecis\'ement, le morphisme ci-dessus est
construit
\`a l'aide du cup-produit avec les puissances de $c_1(\sO_{\P(E)}(1))$; il
suffit de voir que c'est un isomorphisme apr\`es tensorisation par $\Q$
ou par
$\Z/l$ pour un nombre premier $l$ quelconque. Dans le premier cas, c'est
vrai gr\^ace \`a la proposition \ref{p4} par r\'eduction aux groupes de
Chow sup\'erieurs \cite[th. 7.1]{bloch}. Dans le deuxi\`eme cas,  c'est
vrai pour $l\ne \car k$ gr\^ace \`a \cite{geilev} par r\'eduction \`a la
cohomologie \'etale \cite[th. 2.2.1]{jou}, et pour $l=\car k$ gr\^ace \`a
\cite{geilevp} par r\'eduction \`a la cohomologie de Hodge-Witt
logarithmique \cite[cor. I.2.1.12]{gros}.

\subsection{b): images directes} Soit maintenant $f:X\to Y$ un morphisme
projectif entre deux $k$-vari\'et\'es lisses. Nous nous proposons de
construire des morphismes
\[f_*:H^i_\et(X,\Z(n))\to H^{i+2c}_\et(Y,\Z(n+c))\]
pour tout $(i,n)\in\Z\times\Z$, o\`u $c=\dim Y-\dim X$. 

\subsubsection{Cas d'une immersion ferm\'ee}\label{A.7.1} Supposons que
$f$ soit une immersion ferm\'ee. On proc\`ede comme dans
\cite[\S 4]{k}. Supposons d'abord $n\ge 0$. En uti\-li\-sant la
transformation naturelle canonique
$\alpha_X^*f_\Zar^!\to f_\et^!\alpha_P^*$ et le th\'eor\`eme de
localisation pour les groupes de Chow sup\'erieurs, on obtient des
morphismes
\begin{equation}\label{eqA.2}
\Z(n)^\et_X\to Rf_\et^!\Z(n+c)^\et_Y[2c].
\end{equation}

(Ces morphismes ne sont pas des isomorphismes en l'exposant
ca\-rac\-t\'e\-ris\-ti\-que de $k$, \cf \cite[Theorem 4.2 b)]{k}.) Pour
$n$ quelconque, on a par puret\'e de la cohomologie \'etale \`a
coefficients finis des morphismes
\[(\Q/\Z)'(n)_X\iso Rf^!(\Q/\Z)'(n+c)_Y[2c].\]

En composant avec \eqref{eqA.1}, cela fournit des morphismes
\eqref{eqA.2} \'egalement pour $n<0$. Ceci d\'efinit
$f_*$. Par construction, on a 
\begin{equation}\label{eqA.3}
(gf)_*=g_*f_* \text{ si  $f$ et $g$ sont deux
immersions ferm\'ees composables.}
\end{equation}

\'Etant donn\'e un diagramme
cart\'esien
\[\begin{CD}
X'@>f'>> Y'\\
@V{g'}VV @V{g}VV\\
X@>f>> Y
\end{CD}\]
avec $f$ (donc $f'$) immersion ferm\'ee, on a aussi la
\emph{condition d'\'echange}
\begin{equation}\label{eqA.5}
g^*f_*=f'_*(g')^*:H^i_\et(X,\Z(n))\to H^{i+2c}(Y',\Z(n+c)).
\end{equation}

Pour $n<0$, elle r\'esulte de la formule correspondante en cohomologie
\'etale \`a coefficients finis; pour $n>0$ elle s'obtient par
faisceautisation \'etale de la m\^eme formule au niveau des complexes de
faisceaux Zariski.

Par la m\^eme m\'ethode on a l'identit\'e, pour une immersion ferm\'ee
$f$
\[(f\times 1)_*(x\boxtimes y)=f_*x\boxtimes y\]
o\`u $\boxtimes$ d\'esigne le produit externe (``cross-produit"), d'o\`u
l'on d\'eduit\footnote{Je remercie Marc Levine de m'avoir indiqu\'e cette
astuce.} par pull-back via la diagonale la
\emph{formule de projection}
\begin{equation}\label{eqA.4}
f_*(x\cdot f^*y)=f_*x\cdot y.
\end{equation}

\subsubsection{Cas d'une projection} Supposons que $f$ soit la
projection $Y\times \P^{-c}\to Y$. On d\'efinit $f_*$ via la projection
du membre de droite de l'isomorphisme de \ref{fibproj} sur la composante
$j=-c$ du membre de gauche.

\subsubsection{Cas g\'en\'eral} Par hypoth\`ese, $f$ se factorise en
\[X\overset{i}{\To}Y\times \P^n\overset{\pi}{\To} Y\]
o\`u $i$ est une immersion ferm\'ee et $\pi$ est la premi\`ere projection.
On d\'efinit 
$f_*:=\pi_*i_*$.

\subsubsection{Ind\'ependance du choix de la
factorisation}\label{ind} Soit $X\overset{i'}{\To}Y\times
\P^m\overset{\pi'}{\To} Y$ une autre telle factorisation. On a des
immersions ferm\'ees canoniques $a:\P^n\to \P^{m+n+1}$ et $a':\P^m\to
\P^{m+n+1}$. Notons $\pi'':Y\times \P^{m+n+1}\to Y$ la premi\`ere
projection. En tenant compte de
\eqref{eqA.3}, cette construction r\'eduit la question \`a d\'emontrer les
formules $\pi_*=\pi''_*a_*$ et $\pi'_*=\pi''_*a'_*$. En d'autres termes,
on est ramen\'e au cas particulier o\`u $i$ est une inclusion $Y\times
\P^n\to Y\times\P^m$, donn\'ee par l'inclusion canonique $\P^n\to \P^m$
($n\le m$).

Soit $\xi=c_1(\sO_X(1))\in H^2_\et(X,\Z(1))$. L'assertion r\'esulte de la
formule du fibr\'e projectif (appliqu\'ee \`a $X=Y\times \P^n$ et \`a
$Y\times
\P^m$) et de la formule
\[i_*(\xi^j)=
\begin{cases}
0&\text{si $j<-c$}\\
1&\text{si $j=-c$.}
\end{cases}\]

Cette propri\'et\'e est bien connue au niveau des groupes de Chow: elle
r\'esulte par exemple de la fonctorialit\'e covariante de ces derniers
(\cf \cite[prop. 3.1 a)]{fulton}). Mais
$CH^1(X)\iso H^2_\et(X,\Z(1))$, donc elle en r\'esulte pour la
cohomologie motivique \'etale par fonctorialit\'e.

\subsection{(ii): fonctorialit\'e} Soient $f:X\to Y$ et $g:Y\to Z$ deux
morphismes projectifs, munis de d\'ecompositions $f=\pi i$, $g=\pi'i'$
comme ci-dessus. Nous voulons montrer que $(gf)_*=g_*f_*$. Consid\'erons
le diagramme commutatif:
\[\xymatrix{
&Y\times\P^m\ar[d]^\pi\ar[r]^-{i''}& Z\times
\P^m\times\P^n\ar[d]^-{\pi''}\ar[r]^S&
Z\times\P^{mn+m+n}\ar[ddl]^-{\Pi}\\ X\ar[ur]^i\ar[r]^f&
Y\ar[r]^-{i'}\ar[dr]^g& Z\times\P^n\ar[d]^-{\pi'}\\ &&Z
}\]
o\`u le carr\'e central est cart\'esien et $S$ est donn\'e par le
plongement de Segre. Comme par \ref{ind} on a $f_*=\pi_* i_*$,
$g_*=\pi'_*i'_*$ et $(gf)_*=\Pi_* (Si''i)_*$, en tenant de plus compte
de \eqref{eqA.3}, ce diagramme nous r\'eduit
\`a d\'emontrer les deux formules:

\begin{thlist}
\item $i'_*\pi_*=\pi''_* i''_*$;
\item $\pi'_*\pi''_*=\Pi_*S_*$.
\end{thlist}

\subsubsection{Cas de (i)} Soit $\xi=c_1(\sO_{\P^m}(1))$: la formule du
fibr\'e projectif pour $\pi$ et pour $\pi''$ s'exprime en fonction des
pull-backs de $\xi$. L'assertion r\'esulte donc de la formule de
projection \eqref{eqA.4}, appliqu\'ee avec $f=i''$.

\subsubsection{Cas de (ii)} Soient $\xi=c_1(\sO_{\P^m}(1))$,
$\xi'=c_1(\sO_{\P^n}(1))$ et $\xi''=c_1(\sO_{\P^{mn+m+n}}(1))$. En
utilisant la formule du fibr\'e projectif, on a pour $z\in
H^*_\et(Z,\Z(*))$
\[\pi'_*\pi''_*((\pi'\pi'')^*z\cdot \xi^a\cdot (\xi')^b)=
\begin{cases}
0&\text{si $(a,b)\ne (m,n)$}\\
z&\text{si $(a,b)= (m,n)$}
\end{cases}
\]
et de m\^eme
\[\Pi_*(\Pi^*z\cdot (\xi'')^c)=
\begin{cases}
0&\text{si $c\ne mn+m+n$}\\
z&\text{si $c= mn+m+n$.}
\end{cases}
\]

\begin{sloppypar}
En utilisant de nouveau la formule de projection, il s'agit de montrer
que, dans l'\'ecriture de $S_*(\xi^a\cdot (\xi')^b)$ sur la base des
$(\xi'')^c$, le coefficient de $(\xi'')^{mn+m+n}$ est $0$ si $(a,b)\ne
(m,n)$ et $1$ si $(a,b)= (m,n)$. Mais il s'agit l\`a d'un calcul dans
l'anneau de Chow, et le r\'esultat est vrai par fonctorialit\'e
covariante des groupes de Chow.
\end{sloppypar}

\subsection{(v): condition d'\'echange}\label{A.9} Il s'agit d'\'etendre
\eqref{eqA.5} au cas o\`u $f$ (donc $f'$) est projectif. \'Etant donn\'e
\ref{A.7.1}, on se r\'eduit au cas o\`u $f$ est une projection $Y\times
\P^n\to Y$. C'est alors imm\'ediat \'etant donn\'e la formule du fibr\'e
projectif et \ref{A.5}.

\subsection{(vi): formule de projection} Il s'agit de g\'en\'eraliser
\eqref{eqA.4} \`a un morphisme projectif quelconque. M\^eme r\'eduction
et justifications qu'en \ref{A.9}.

\subsection{Cas de la cohomologie de Lichtenbaum
(justification de la proposition \protect\ref{p5.8})}\label{A.11} M\^emes
m\'ethodes, en \'etendant la formule du fibr\'e projectif \`a cette
cohomologie. Pour ceci, on construit des morphismes comme dans
\ref{fibproj} pour la cohomologie de Lichtenbaum. On d\'emontre que ce
sont des isomorphismes par r\'eduction au cas \'etale gr\^ace
\`a \eqref{e1}.

\end{document}